\documentclass[a4paper]{amsart}
\usepackage[a4paper, margin={3cm}]{geometry}

\usepackage{graphicx}
\usepackage{epstopdf} 
\usepackage{bm}
\usepackage{todonotes}
\usepackage{booktabs}
\usepackage{caption}
\usepackage{subcaption}
\usepackage{placeins}
\usepackage{comment}
\usepackage{cancel}
\usepackage{multirow}
\usepackage{pdflscape}
\graphicspath{{figures/}}
\usepackage{hyperref}

\newcommand{\rd}{\mathrm{d}}
\newcommand{\vt}[1]{\bm{#1}}

\newcommand{\dd}[2]{\frac{\partial #1}{\partial #2}}
\newcommand{\q}{\vt{u}} 
\newcommand{\qbar}{\bar{\q}} 
\newcommand{\pphys}{\vt{\mu}} 
\newcommand{\pcal}{\vt{\theta}} 
\newcommand{\G}{\vt{G}} 
\newcommand{\qred}{\bar{\vt{v}}} 
\newcommand{\F}{\vt{F}} 
\newcommand{\f}{\vt{f}} 
\newcommand{\Fbar}{\bar{\F}} 
\newcommand{\m}{\vt{m}} 
\newcommand{\mdiv}{\hat{\vt{m}}} 
\newcommand{\me}{\vt{m}^{e}} 
\newcommand{\ms}{\vt{m}^{s}} 
\newcommand{\C}{\vt{C}} 
\newcommand{\Jof}{J^{\mathrm{prio}}} 
\newcommand{\Jtf}{J^{\mathrm{post}}} 
\newcommand{\Ltf}{\mathcal{L}^{\mathrm{post}}} 
\newcommand{\Lof}{\mathcal{L}^{\mathrm{prio}}} 
\newcommand{\Jr}{J^{\mathrm{r}}} 
\newcommand{\R}{\mathcal{R}}

\newcommand{\Fh}{\vt{F}_{h}} 
\newcommand{\qh}{\vt{u}_{h}} 
\newcommand{\Gh}{\vt{G}_{h,\pcal}} 
\newcommand{\qhred}{\bar{\vt{v}}_{h}} 
\newcommand{\fh}{\vt{f}_{h}} 
\newcommand{\mh}{\vt{m}_{h,\pcal}} 
\newcommand{\meh}{\vt{m}^{e}_{h,\pcal}} 

\newcommand{\obs}{\vt{y}} 

\newcommand{\rh}{\vt{r}_{h}}

\newcommand{\ro}{\vt{r}}

\newcommand{\pd}[2]{\frac{\partial#1}{\partial#2}}

\DeclareMathOperator*{\argmin}{arg\,min}
\newcommand{\RC}[1]{\label{#1}} 
\newcommand{\revone}[1]{\textcolor{black}{#1}}
\newcommand{\revtwo}[1]{\textcolor{black}{#1}}
\newcommand{\revthree}[1]{\textcolor{black}{#1}}

\begin{document}

\title{Scientific machine learning for closure models in multiscale problems: a review} 

\author{B. Sanderse}
\address{Scientific Computing Group, Centrum Wiskunde \& Informatica, Science Park 123, Amsterdam, The Netherlands}
\email{b.sanderse@cwi.nl}
\author{P. Stinis}
\address{Advanced Computing, Mathematics and Data Division, Pacific Northwest National Laboratory, Richland, 99352, WA, USA}
\email{panagiotis.stinis@pnnl.gov}
\author{R. Maulik}
\address{College of Information Sciences and Technology, The Pennsylvania State University, Westgate Building, University Park, 16802, PA, USA}
\email{rmaulik@psu.edu}
\author{S. E. Ahmed}
\address{Advanced Computing, Mathematics and Data Division, Pacific Northwest National Laboratory, Richland, 99352, WA, USA}
\email{shady.ahmed@okstate.edu}

\begin{abstract}
Closure problems are omnipresent when simulating multiscale systems, where some quantities and processes cannot be fully prescribed despite their effects on the simulation's accuracy.
Recently, scientific machine learning approaches have been proposed as a way to tackle the closure problem, combining traditional (physics-based) modeling with data-driven (machine-learned) techniques, typically through enriching differential equations with neural networks. This paper reviews the different reduced model forms, distinguished by the degree to which they include known physics, and the different objectives of a priori and a posteriori learning. The importance of adhering to physical laws (such as symmetries and conservation laws) in choosing the reduced model form and choosing the learning method is discussed. The effect of spatial and temporal discretization and recent trends toward discretization-invariant models are reviewed. In addition, we make the connections between closure problems and several other research disciplines: inverse problems, Mori-Zwanzig theory, and multi-fidelity methods. In conclusion, much progress has been made with scientific machine learning approaches for solving closure problems, but many challenges remain. In particular, the generalizability and interpretability of learned models is a major issue that needs to be addressed further.
\end{abstract}

\maketitle


\section{Introduction and background}
Closure problems are everywhere. They arise when a mathematical-physical model that describes certain quantities of interest is created, \textit{which depends on quantities that are not of prime interest but whose effect cannot be neglected}. A prominent example is in numerical weather prediction, where the aim is to predict tomorrow's weather given the physical laws that govern fluid flows and today's conditions. In principle, these physical laws describe all relevant physics, but in practice certain processes like cloud formation take place on such small spatial and temporal scales that resolving (simulating) them on a computer is impossible, even though \textit{their effect on the weather is crucial} \cite{christensen2022}. 

The weather prediction problem is an example of a multiscale problem in which the quantities of interest are associated with large scales, but depend on quantities associated with small scales. If these small scales are not resolved, the large-scale equations are \textit{unclosed}. 
In this article, we focus on multiscale problems that lead to a closure problem, and in particular we consider multiscale problems that involve a continuum of scales (fluid flows being an important example). There is no clear separation of scales, and the closure model is needed everywhere in the domain to supplement the large-scale model (Type ``B'' problems in the terminology of E \cite{e2011}).

Closure models appear under different names in different domains \cite{ahmed2021closures}. The term ``closure'' is most common in the fluid dynamics community, designating the turbulence closure problem---closure models are known as subgrid or subgrid-scale (SGS) models. In kinetic equations and radiation transport, the term ``moment closure'' is used \cite{levermore1996,han2019,huang2022}. In numerical weather prediction and climate science, the problem is known as ``parameterization'', and finding accurate and stable parameterizations is considered a key challenge \cite{christensen2022,frezat2023}. In general circulation models (including atmospheric and oceanic models), the same terminology is used \cite{sonnewald2021,bolton2019}. In the reservoir simulation community, the closure problem is commonly referred to as ``homogenization'' or ``upscaling'' \cite{christie1996}. In materials science, molecular dynamics, and computational biology, the term ``coarse graining'' is used \cite{saunders2013,giessen2020}.

Finding a closure model for multiscale problems is difficult and often more an art than a science. As mentioned, part of the model description is often known from physical laws (such as conservation laws), and only the effect of the small scales on the large scales needs to be approximated. The mapping associated with this approximation is typically strongly nonlinear---therefore, machine learning methods form a promising avenue, given their strong approximation properties. Consequently, a hybrid approach to the closure problem has become a popular research field, involving a combination of a physics-based model describing the large scales and a machine-learning based model to approximate the effect of the small scales. Such approaches fall within the realm of what is now known as scientific machine learning. Because the physics-based model typically consists of differential equations, this hybrid approach in practice means a mixture of differential equations with neural networks. This leads to many new and interesting research questions regarding stability, consistency, and convergence, which have only started to be addressed recently. 

Consequently, the aim of this review is to give an overview of recent scientific machine learning approaches to the closure problem and the open challenges faced by the community. To do so, we give a somewhat abstract view on the closure problem, which is useful because it provides a generic description that fits several application domains. A second aim of this review is to make the connection to adjacent research disciplines clear; for example, the inverse problem community, the Mori-Zwanzig community, etc. We realize that our notation and the cited literature are somewhat geared toward fluid flow applications (and large eddy simulations [LESs] in particular), but we have tried to make it general enough to address a broad audience. We also note that many excellent review articles have appeared on the topic of machine learning, partial differential equations and/or fluid mechanics, for example \cite{duraisamy2019,ahmed2021closures,brunton2023,lino2023,viquerat2022}. In this review, we therefore assume that the reader is familiar with basic concepts of machine learning and differential equations.


The outline of this paper is as follows. Section \ref{sec:problem_formulation} and Section \ref{sec:objective_functions} describe the fundamentals of closure modeling, splitting the closure problem into two parts: the reduced model form (\textit{what} to learn), and the choice for the objective function (\textit{how} to learn it). Together with Section \ref{sec:physics-constrained} (embedding physics constraints) and Section \ref{sec:discretization} (spatial and temporal discretization aspects), these sections describe the methodological framework in which many existing closure model approaches fit. Sections \ref{sec:inverse_etc}, \ref{sec:MZ} and \ref{sec:multifidelity} each give a different perspective on the closure problem from the viewpoint of inverse problems, from the viewpoint of Mori-Zwanzig theory, and from the viewpoint of multi-fidelity and multiscale methods. Section \ref{sec:conclusions} concludes with directions for future research.


\section{Reduced model forms}\label{sec:problem_formulation}

\subsection{Model reduction}
We will consider partial differential equations (PDEs) of the form
\begin{equation}\label{eqn:PDE}
\F(\q;\pphys) = 0,
\end{equation}
where $\F$ is the PDE operator (containing temporal and spatial derivatives), $\q(\vt{x},t)$ is the solution, and $\pphys$ represents physical parameters or parameterized initial or boundary conditions. 
Equation \eqref{eqn:PDE} will be referred to as the \textit{full model}. In practice, this full model is solved by discretizing in space and time. The resulting discretized model will be referred to as the \textit{high-fidelity model}, \RC{rev3_comment1}\revthree{assuming a sufficient resolution is maintained to capture the relevant spatial and temporal scales.}. 

The degrees of freedom or spectral content of the solution $\q(\vt{x},t)$ can be \textit{reduced}, e.g.,\ through filtering, averaging, projection, or truncation. The reduction operator $\mathcal{A}$ leads to a reduced field $\qbar(\vt{x},t)$: 
\begin{equation}\label{eqn:qbar_def}
\qbar := \mathcal{A} (\q).
\end{equation}
$\mathcal{A}$ is often a linear operator, but this is not necessary. We assume in this review that \textit{the main interest is the computation of an accurate approximation to $\qbar$}, not in computing an accurate approximation to $\q$. We also assume that a small number of simulations of \eqref{eqn:PDE} can be performed, leading to reference data for both the solution field and the reduced field. Note that a further reduction is sometimes made, i.e.,\ in certain applications where one is only interested in a quantity of interest $Q(\qbar)$ that is of lower dimension than $\qbar$ itself. Two straightforward but important relations in the context of model reduction are \revthree{$\Fbar(\q;\pphys) := \mathcal{A}(\F(\q;\pphys))=0$} \footnote{assuming \revthree{$\mathcal{A}$ is linear}} and $\F(\qbar;\pphys)\neq 0$. In other words, the reduction operator $\mathcal{A}$ and the PDE operator $\F$ do not necessarily commute.

The reduction operator $\mathcal{A}$ can be local in space and time, but also nonlocal (e.g.,\ to account for memory effects). A prominent example of a reduction operator is a filter with a kernel $K$ \cite{sagaut2006}:
\begin{equation}\label{eqn:continuous_filter}
\qbar(\vt{x},t) = \int \q(\vt{\xi},t) K(\vt{x},\vt{\xi}) \rd \vt{\xi},
\end{equation}
where the kernel is often of the convolution type, so $K(\vt{x},\vt{\xi})=K(\vt{x}-\vt{\xi})$. 
Another example of a reduction operator is an expansion in an orthogonal basis $\vt{\phi}_{i}(\vt{x})$ truncated to a finite dimension:
\begin{equation}\label{eqn:continuous_projection}
    \qbar(\vt{x},t) = \sum_{i=1}^{N} a_{i}(t) \vt{\phi}_{i}(\vt{x}), 
\end{equation}
which is commonly employed in the Reduced-Order Model (ROM) community \cite{hesthaven2016,ahmed2021closures}.

The main problem is that $\qbar$ is in general not a solution to PDE \eqref{eqn:PDE}. The common approach to predict $\qbar$ is to formulate a parameterized \textit{reduced model} $\G$ (typically a PDE, but this is not necessary) with parameters $\pcal$, which describes an approximation $\qred \approx \qbar$ \textit{that does not \revthree{explicitly} depend on $\q$}:
\begin{equation}\label{eqn:filteredPDE}
    \boxed{\G_{\pcal}(\qred;\pphys) = 0}.
\end{equation}
This model $\G$ should be computationally much faster to solve than the original PDE \eqref{eqn:PDE}, while retaining most of the accuracy \RC{rev3_comment6}\revthree{of the quantity of interest (QoI)}. We switched notation from $\qbar(\vt{x},t)$ to $\qred(\vt{x},t)$ because solving \eqref{eqn:filteredPDE} will in general not yield the reduced field $\qbar$, since $\qbar$ requires knowledge of the full solution $\q$, as defined by \eqref{eqn:qbar_def}. In addition, $\qred(\vt{x},t)$ could possibly include additional degrees of freedom, such as latent variables, that aid in approximating the evolution of $\qbar$. The vector $\pcal$ contains parameters or parameterized functions (e.g.,\ neural networks, polynomials, etc.) that need to be \RC{rev3_comment7}\revthree{learned or discovered, depending on the hypotheses and Ansatz presumed in the problem formulation.} 


The main challenge lies in deriving an expression for $\G$ (including the specification of all necessary parameters), thus placing us in the setting of \textit{model discovery}  \cite{brunton2016discovering,rudy2017,loiseau2018sparse,vaddireddy2020feature,fukami2021sparse} \RC{rev2_comment1}\revtwo{(even though the term model discovery is often used in the context of discovering the PDE for $\q$ instead of for $\qbar$)}. This means that the governing equations are not fully known, and techniques such as physics-informed neural networks (PINNs) \cite{raissi2019} and neural finite element methods \cite{e2018} cannot be easily applied because these assume that the PDE is known (it enters in the loss function). For a review of the development of these methods, see for example\ \cite{lino2023,kumar2011,cai2021,karniadakis2021physics}. However, because PINNs can also be used to solve inverse problems, they might be effective in detecting closure models if the closure model is framed as an inverse problem---see sections \ref{sec:soft_constraints} and \ref{sec:inverse}. 

The difficulty of deriving an expression for $\G$ is not only due to its inverse-problem nature, but also stems from the fact that the dynamics described by the full model $\F$ is often nonlinear and \textit{chaotic}. This means that the full model $\F$ would typically be solved with a statistical approach, and a similar strategy can be expected for the reduced model---see Sections \ref{sec:statistical} and \ref{sec:MZ}. Nevertheless, many reduced models used to date are still deterministic. 

\subsection{Closure problem}
Equation \eqref{eqn:filteredPDE} describes a very generic class of models for $\qred$. A natural approach is to find $\G$ such that it approximates the reduced model $\Fbar$, i.e.,\ $\G \approx \Fbar$. Applying the reduction operation to \eqref{eqn:PDE} yields:
\begin{equation}
    \Fbar(\q;\pphys):=\mathcal{A}(\F(\q;\pphys)) = 0 \quad \rightarrow \quad \F(\qbar;\pphys) +  \C(\q,\qbar;\pphys) = 0
\end{equation}
where 
\begin{equation}\label{eqn:commutation_error}
    \C(\q,\qbar;\pphys) = \Fbar(\q;\pphys) - \F(\qbar;\pphys). 
\end{equation}
This equation shows that one cannot simply solve the original PDE for the reduced field, because the reduction step and the PDE operator do not commute. Generally, the equation for $\Fbar$ involves $\q$ and is therefore \textit{unclosed}. The common expression for $\G$ follows from approximating the commutator error $\C$ by a model $\m_{\pcal}(\qbar;\pphys)$ that depends only on $\qbar$:
\begin{equation}
    \m_{\pcal}(\qbar;\pphys) \approx \C(\q,\qbar;\pphys),
\end{equation}
so that
\begin{equation}\label{eqn:closure_model_form}
    \boxed{\G_{\pcal}(\qred;\pphys) = \F(\qred;\pphys) +  \m_{\pcal}(\qred;\pphys)} .
\end{equation}
The model $\m_{\pcal}$ is known as the \textit{closure model}. Expression \eqref{eqn:closure_model_form} is the basic equation for the reduced model, which expresses that the closure model approximates a commutator error that acts as an additive correction to the full model.  If additional knowledge of $\F$ is available, e.g.,\ when $\F$ is of the form 
\begin{equation}\label{eqn:evolution_form}
    \F(\q;\pphys) = \dd{\q}{t} + \f(\q;\pphys),
\end{equation}
the expression for $\vt{G}$ reads \footnote{assuming the reduction operator commutes with temporal differentiation}:
\begin{equation}\label{eqn:closure_model_form_ddt}
    \boxed{\G_{\pcal}(\qred;\pphys) = \dd{\qred}{t} + \f(\qred;\pphys) + \m_{\pcal}(\qred;\pphys)} . 
\end{equation}
In practice, this form is as important as equation \eqref{eqn:closure_model_form} because in many physical problems, $\F$ is of the form \eqref{eqn:evolution_form}. In addition, in fluid flow and material science applications, $\f$ usually involves the divergence of a stress tensor, and \eqref{eqn:closure_model_form_ddt} is further specialized to $\G_{\pcal}(\qred;\pphys) = \dd{\qred}{t} + \f(\qred;\pphys) + \nabla \cdot \mdiv_{\pcal}(\qred;\pphys)$. An example is given in Section \ref{sec:LES}. 

In \RC{rev3_comment10a}\revthree{many} scientific machine learning approaches, the term $\m$ is modeled by a neural network. This is a logical choice, given that the mapping from $\qbar$ to $\m(\qbar)$ is generally highly nonlinear. The resulting PDE consists of ``known physics'', $\f(\qred)$, plus ``learned physics'' (a neural network). However, the stability of this PDE solution is not guaranteed and is a subject of significant research, which will be further discussed in Section \ref{sec:objective_functions}. \RC{rev3_comment10b} \revthree{In addition, approaches like reinforcement learning provide an alternative perspective on learning $\m$ -- see section \ref{sec:RL}.}

\subsection{Other reduced model forms}\label{sec:other_reduced_forms}
If one removes the condition that $\G$ should approximate $\F$ and the focus is only on having $\qred \approx \qbar$, many other forms for $\G$ are possible, typically depending on how much a priori knowledge of $\F$ is given. In addition, the degree of code access (e.g.,\ whether the source code is available, or only a compiled version) can determine the form chosen for $\G$. These non-closure forms are not the main focus of this paper but will be mentioned for completeness. For example, if the exact form of $\F$ is not known or available and one only knows that the reduced model should be of evolution form, one can learn the \textit{entire} right-hand side (instead of only learning an approximation to the commutator error):
\begin{equation}\label{eqn:neuralODE}
    \G_{\pcal} (\qred;\pphys) = \dd{\qred}{t} + \me_{\pcal}(\qred;\pphys) .
\end{equation}
The symbol $\me$ is used to distinguish it from $\m$ because $\me$ generally does not approximate the commutator error. Form \eqref{eqn:neuralODE} is very generic, and many equation-discovery methods for dynamical systems fit in this form \cite{rudy2017}. A recent example for turbulent flow simulation is \cite{stachenfeld2022a}, where the approach is called ``fully learned''---no known physics are used in choosing the form of the reduced model, except that the problem is a dynamical system. Recent weather prediction models employ a similar setup \cite{pathak2022,lam2023}. If \RC{rev3_comment12a}\revthree{$\me_{\pcal}$} is represented by a neural network, this form is similar to neural ODEs, and this will be discussed in Section \ref{sec:neuralODE}. 

If even less a priori knowledge is available, one can resort to so-called \textit{surrogate models}, 
\begin{equation}\label{eqn:G_surrogate}
    \G_{\pcal} (\qred;\pphys) = \qred(\vt{x},t) + \ms_{\pcal}(\vt{x},t;\pphys) , 
\end{equation}
where \RC{rev3_comment12b}\revthree{$\ms_{\pcal}$} represents a linear or nonlinear function approximator, such as neural networks, neural operators, polynomial expansions, Gaussian processes, \RC{rev3_comment13}\revthree{dynamic mode decomposition, and proper orthogonal decomposition,} taking physical parameters and spatial and temporal coordinates as input instead of the solution $\qred$. Examples of neural operators that fit in this class are the Deep Operator Network (DeepONet) \cite{lu2021}, Fourier Neural Operator (FNO) \cite{li2021}, Laplace Neural Operator (LNO) \cite{cao2023lno}, Graph Neural Operator (GNO) \cite{xu2024equivariant}, and Convolutional Neural Operator (CNO) \cite{raonic2024}.

Next to the closure form \eqref{eqn:closure_model_form}, the generic evolution form \eqref{eqn:neuralODE}, and the surrogate form \eqref{eqn:G_surrogate}, many other reduced forms are  possible. For example, equation \eqref{eqn:closure_model_form_ddt} for the reduced variable can be extended with equations for latent variables, which can be used to model the effect of memory (for details, see Section \ref{sec:MZ}). In the context of learning ODEs, this leads to so-called augmented neural ODEs \cite{dupont2019}, which have been used in the context of closure models \cite{melchers2022,vangastelen2023}. \revthree{In the context of the equation-free methods \cite{kevrekidis2003}, the learned effective dynamics (LED) framework has been proposed for learning low-dimensional models that include memory \cite{vlachas2022}.} Another example of a reduced form is a stochastic differential equation with trainable terms \cite{tzen2019,boral2023}.


\subsection{Example from fluid dynamics: large eddy simulation}\label{sec:LES}
An important example of closure models appears in the simulation of turbulent flows with LES \cite{sagaut2006,berselli2006a}. In this case, the incompressible Navier-Stokes equations form the high-fidelity model,
\begin{equation}
    \F(\q;\pphys) := \dd{\q}{t} + \nabla \cdot (\q \otimes \q) - \nabla \cdot (2\nu \vt{S}(\q)) + \nabla p,
\end{equation}
where $\vt{u}$ is the velocity, $p$ the pressure, and $\vt{S}(\q)=\tfrac{1}{2}(\nabla \q + (\nabla \q)^T)$ is the strain-rate tensor. \RC{rev3_comment15} \revthree{The term $\nabla \cdot (\q \otimes \q)$ represents non-linear convection (with $\q \otimes \q$ the dyadic product) and is the major cause of the closure problem. It is often written in the advective form $(\q \cdot \nabla) \q$.} The equation is supplemented with the constraint $\nabla \cdot \q = 0$, but to keep the discussion concise, we omit this constraint in what follows. The reduction operator is a convolutional filter, similar to \eqref{eqn:continuous_filter}, which commutes with spatial and temporal differentiation so that the expression for the commutator error reads
\begin{equation}\label{eqn:commutator_LES}
    \C(\q,\qbar) = \Fbar(\q;\pphys) - \F(\qbar;\pphys) = \nabla \cdot (\overline{\q \otimes \q} - (\qbar \otimes \qbar)).
\end{equation}
Because the commutator error is the divergence of a quantity, the aim is to construct an approximation $\mdiv$ such that
\begin{equation}\label{eqn:commutator_tensor}
    \C(\q,\qbar) \approx \nabla \cdot \mdiv_{\pcal} (\qbar;\pphys).
\end{equation}
A classical model for $\mdiv$ is the Smagorinsky model,
\begin{equation}\label{eqn:smagorinsky}
    \mdiv_{\pcal}(\qbar;\pphys) =  - \nu_{t}(\theta) \vt{S}(\qbar),           
\end{equation}
and the goal is to learn the Smagorinsky coefficient $\theta=C_S$, where $\nu_{t}(\theta) = (\theta \Delta)^2 |\bar{\vt{S}}|$ and $\Delta$ denotes the filter width \cite{berselli2006a}. The reduced model expression is in this case
\begin{equation}\label{eqn:closure_model_form_Smagorinsky}
    \G_{\pcal}(\qred;\pphys) := \dd{\qred}{t} + \nabla \cdot (\qred \otimes \qred) - \nabla \cdot (2 \nu \vt{S}(\qred)) + \nabla \bar{p} - \nabla \cdot( \nu_{t}(\theta) \vt{S}(\qred)).
\end{equation}
Note that even though the Reynolds-averaged Navier-Stokes equations (RANS) have a similar form, they feature a different reduction operator \RC{rev3_comment17} \revthree{that leads to a much larger commutator error, and it should not be confused with LES}. In many scientific machine learning approaches, instead of learning the single parameter $C_{S}$, a neural network is trained to approximate $\mdiv$, or $\mdiv$ is expressed in terms of a tensor basis with coefficients that need to be learned (more details in Section \ref{sec:physics-constrained}). 

\section{Objective function choice: a priori vs.\ a posteriori learning}\label{sec:objective_functions}
Having introduced different approaches to construct a reduced model, \RC{rev3_comment18}\revthree{it is clear that the resulting form involves a parameterized model with unknown parameters. This leads to the formulation of an inverse problem, where one or more parameters are learned to meet certain requirements.} In this section, we continue with introducing three distinct methods for approaching the closure problem in terms of the objective function that is being minimized.

\subsection{Approach 1: A priori learning}\label{sec:operator_fitting}
The most common approach in closure modeling is to find the parameters $\theta$ by minimizing an error metric that can be evaluated ``offline'', i.e.,\ \textit{without solving the reduced model} \eqref{eqn:filteredPDE}. This approach is known as \textbf{a priori learning}, ``direct approach'' \cite{macart2021}, or ``offline mode'' \cite{rasp2020}. The optimization problem has the following form:
\begin{equation}\label{eqn:closure_problem_of}
\pcal^{*} = \argmin_{\pcal} \Jof_{\pcal}(\qbar,\q),
\end{equation}
In this approach, a limited number of reference trajectories (solutions of $\q$ as a function of time for different parameters $\pphys$) are obtained from  \eqref{eqn:PDE}, from which $\qbar$ can be extracted  by applying the filter \eqref{eqn:continuous_filter}. 
The objective function typically involves the summation over a set of training data $\mathcal{T}$ \RC{rev3_comment19a}\revthree{(in parameter space)},
\begin{equation}\label{eqn:loss_function_sum}
    \Jof_{\pcal} (\qbar,\q) = \frac{1}{|\mathcal{T}|} \sum_{i \in \mathcal{T}} \Lof_{\pcal}(\qbar(\vt{x},t;\pphys_{i}),\q(\vt{x},t;\pphys_{i})),
\end{equation}
where $|\mathcal{T}|$ measures the cardinality (``size'') of the training data set, \RC{rev3_comment19b}\revthree{and $\pphys_{i}$ are the parameters corresponding to the $i$-th training sample, e.g.\ a vector describing initial conditions. The summation can be easily extended to include a temporal aspect}. The loss function $\Lof$ can have many different forms. A common example is to minimize the squared commutator error (see closure model form \eqref{eqn:closure_model_form}):
\begin{equation}\label{eqn:Lot_2}
      \Lof_{\pcal}(\qbar,\q)  = \| \m_{\pcal} (\qbar;\pphys) - \C(\q,\qbar;\pphys) \|^2,
\end{equation}
\RC{rev3_comment20}\revthree{where the \emph{true} commutator error $\C$ is computed according to \eqref{eqn:commutation_error} given the collected data set of $\q$ and its filtered value $\qbar$.} Because in this case one learns an operator, this approach is sometimes called operator fitting, or operator inference in the ROM community \cite{peherstorfer2016,kramer2024} (note that this term should not be confused with the neural operators mentioned in Section \ref{sec:other_reduced_forms}). An example of such an operator is the commutator error in LES, equation \eqref{eqn:commutator_LES}. The general expression for the loss function is
\begin{equation}\label{eqn:Lot_1}
     \boxed{\Lof_{\pcal}(\qbar,\q) = \| \G_{\pcal} (\qbar;\pphys) - \Fbar(\q;\pphys)) \|^2}.
\end{equation}
Figure \ref{fig:a-priori_a-posteriori} (left) shows how in a priori learning, \textit{solving} the reduced model $\G=0$ is not required in the optimization process; one only needs to \textit{evaluate} the residual $\G_{\pcal}(\qbar;\pphys)$. 
\RC{rev3_comment21}\revthree{Since $\Fbar=0$, expression \eqref{eqn:Lot_1} effectively constitutes the residual of the approximate model $\G$, evaluated using the reference trajectory $\qbar$). This means that a priori learning amounts to \textit{residual minimization}. To evaluate the residual, the full state $\q$ is not required, but only the reduced state $\qbar$.}\\ \revthree{Note that keeping $\Fbar(\q)$ in expression \eqref{eqn:Lot_1} is insightful. Consider, for example, $\F$ and $\G$ in the form of equations \eqref{eqn:evolution_form} and \eqref{eqn:closure_model_form_ddt}
---the residual then reads $\G_{\pcal} (\qbar;\pphys) - \Fbar(\q;\pphys) = \m_{\pcal}(\qbar;\pphys) - (\bar{\f}(\q;\pphys)- \f(\qbar;\pphys))$, which shows that $\m$ models a commutator error.}



A priori learning with operator fitting is a common approach in the closure modeling community, and many examples are available: learning stress tensors (or divergence of stress tensor) in fluid dynamics \cite{gamahara2017,ling2016,ling2016a,xiao2016,wang2017,schmelzer2018,beck2019,edeling2015,yang2019a,park2021,maulik2019,kurz2022,guan2023}, learning a polymer stress tensor \RC{rev2_comment2}\revtwo{or} the potential function in molecular dynamics \cite{e2023}, learning closure terms by representing them using a neural operator \cite{lutjens2022}, and parameterization in climate models \cite{rasp2018}.

The crucial aspect of a priori learning is that evaluating the objective function \eqref{eqn:closure_problem_of} only requires reference data $\q$ and $\qbar$, \RC{rev3_comment22} \revthree{obtained as solution vector from solving the high-fidelity model}. It involves the accuracy of the operator, but not the accuracy of the solution $\qred$ that will be computed when the learned operator is employed (which is known as the ``online'' mode \cite{rasp2020}, or the ``a posteriori'' analysis). A common problem is that even with a highly accurate operator fit, the solution $\qred$ can \textit{drift} from the true solution $\qbar$ or becomes \textit{unstable}. The issue is also known under the term \textit{model-data inconsistency} \cite{duraisamy2021,kurz2021}, where the data used to train the closure models is not consistent with the environment in which the model will be run, e.g., due to discretization effects. 

In turbulence, the instability is often associated with the concept of backscatter: energy transfer from the small scales back to the large scales \cite{guan2022b}. For example, \cite{park2021} reported unstable results and applied clipping to limit backscatter. Beck et al.\ \cite{beck2019} performed a projection of their neural-network closure model onto an eddy viscosity basis to enforce stability, and Kurz et al.\ \cite{kurz2021} performed ``stability training'' by adding noise to the training data. \RC{rev2_comment4}\revtwo{This approach of adding noise is common in the machine learning communities, see for example \cite{stachenfeld2022a} and \cite{pfaff2021}.} Rasp et al.\ \cite{rasp2020} proposed coupled online learning in which the (offline-trained) machine learning model is further corrected by nudging with a high-fidelity model that is run in parallel. Pawar and San used data assimilation techniques to correct the trained model in an online environment \cite{pawar2021}. Yuval et al.\ \cite{yuval2020} address the stability issue by employing random forests, which automatically respect energy conservation. Charalampopoulos et al.\ \cite{charalampopoulos2022} used imitation learning to stabilize their results. Guan et al.\ \cite{guan2022b} showed that increasing the size of the training set can give stable results. Jakhar et al.\ \cite{jakhar2023} used a sparsity-promoting linear regression technique to construct interpretable closed-form equations. \RC{rev3_comment23a}\revthree{A more general example of using a-prior closures on dynamical systems was presented in \cite{chekroun2024high} where it was observed that such techniques did not reproduce extreme events from the original system.} Pedersen et al.\ \cite{pedersen2023} proposed to add a term into the a priori loss function that measures the mismatch between the true reduced model state and a prediction by another neural network (called an emulator). By using this emulator instead of the actual reduced model, it is possible to include the solution error into the loss function without needing to resort to a posteriori learning. \RC{rev3_comment23b}\revthree{Similar ideas (i.e., using high-fidelity data to inform parameterizations in low-order models) are also observed in other multiscale modeling applications such as for the mechanical analysis of materials \cite{kirchdoerfer2016data,karapiperis2021data}.}




In summary, the main advantages of the a priori approach are the relative ease of training (no differentiable solvers needed) and the prospect of generalization and interpretability, i.e.,\ the possibility to reuse the learned operator when geometry, boundary conditions, or even the differential equation itself changes. However, the main disadvantage is that the solution of the reduced model is not part of the error metric, so that instability and drift can lead to inaccurate solutions. \RC{rev3_comment24} \revthree{We note that in the reduced order modeling literature, conditions have been derived under which non-intrusive models recover the same models that would be obtained with intrusive projection-based model reduction \cite{kramer2024}. Similar techniques might be useful to derive conditions under which the learned operators in a-priori learning can be reused.}


\begingroup
\begin{figure}[hbtp]
\fontfamily{lmss} 
\fontsize{9pt}{9pt}\selectfont
\centering 
\def\svgwidth{\textwidth}
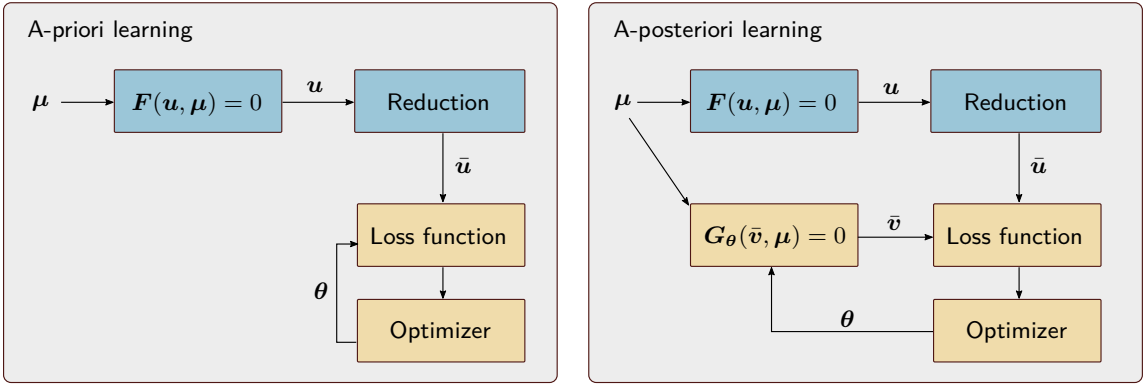 
\caption{A priori vs.\ a posteriori learning. \revthree{Training data $\q$ and $\qbar$ are obtained from solving the high-fidelity model $\F$.} In blue: simulation of the high-fidelity model, from which the ground truth $\qbar$ is derived. In yellow: training. In a posteriori learning, the loss function \eqref{eqn:Ltf} includes evaluating the gradients of the reduced model, while in a priori learning, only the residual must be evaluated (equation \eqref{eqn:Lot_1}).
\label{fig:a-priori_a-posteriori}}
\end{figure}
\endgroup

\subsection{Approach 2: A posteriori learning}\label{sec:solution_fitting}
An alternative approach that addresses some of the issues of a priori learning is to minimize an error metric that involves querying the model being learned itself to compute $\qred$. This amounts to solving the following minimization problem
\begin{equation}\label{eqn:closure_problem_tf}
\pcal^{*} = \argmin_{\pcal} \, \Jtf_{\pcal} (\qred,\qbar),
\end{equation}
where $\Jtf$ is a suitable objective function, and $\qred$ satisfies the reduced model equation \eqref{eqn:filteredPDE}. 
If the reduced model equation is a PDE, then this is effectively a PDE-constrained optimization problem or an optimal control formulation, as will be detailed in Section \ref{sec:optimalcontrol}. Compared to a priori learning, the difference lies in the arguments: $\qbar$ and $\qred$, instead of $\q$ and $\qbar$ in \eqref{eqn:closure_problem_of}. This approach is called \textbf{a posteriori learning} because it involves solving the reduced model equation---see Figure \ref{fig:a-priori_a-posteriori} (right).

The objective function $\Jtf$ again involves summation of a loss function over a set of training data. In addition to the possibility of residual minimization \eqref{eqn:Lot_1}, a posteriori learning offers a second possibility, namely solution error minimization:
\begin{equation}\label{eqn:Ltf}
\Ltf_{\pcal} (\qred,\qbar)= \begin{cases} 
\| \G_{\pcal} (\qred;\pphys) - \Fbar(\q;\pphys)) \|^2, & \text{residual minimization},\\
\| \qred_{\pcal} - \qbar \|^2, & \text{solution error minimization}. 
\end{cases} 
\end{equation}
\RC{rev3_comment26}\revthree{These two loss functions can be used individually, but can also be combined into a single loss function with a weighting factor.} An example of residual minimization is again the commutator loss $\| \m_{\pcal} (\qred;\pphys) - \C(\q,\qbar;\pphys) \|^2$, but now with $\G$ depending on $\qred$, not $\qbar$ as in \eqref{eqn:Lot_1}. An example of solution error minimization is the integral loss $\int_{0}^{T} \| \qred_{\pcal}(\cdot,t) - \qbar(\cdot,t) \|^2  \rd t$, which is sometimes known as ``trajectory fitting'' \cite{melchers2023,crilly2024,chen2022finding}. \RC{rev2_comment3} \revtwo{We note that pure dependence on $\qbar$, and not on $\q$, as suggested in \eqref{eqn:closure_problem_tf} is not always feasible, and the residual minimization in \eqref{eqn:Ltf} is an example of this.}


Again, compared to a priori learning \eqref{eqn:closure_problem_of}, $\Ltf$ in \eqref{eqn:Ltf} is also a function of $\qred$ (which depends on $\pcal$), whereas $\Lof$ was only a function of the reference data $\q$ and $\qbar$, \textit{which do not depend on} $\pcal$. This makes the optimization problem more difficult to solve, and differentiable solvers are typically needed. For example, in order to compute $\dd{J}{\pcal}$, the Jacobian $\dd{\qred}{\pcal}$ is needed---this requires that the code that solves $\G_{\pcal}=0$ is differentiable, or that an adjoint solver is available. In addition to a posteriori learning, the approach is known under various other names, such as ``solver-in-the-loop'' \cite{list2022a}, ``curriculum training'' \cite{keisler2022forecasting}, ``indirect approach'' \cite{macart2021}, ``end-to-end learning'', ``differentiable physics'', and ``online learning'' \cite{frezat2023}. 

The a posteriori learning framework is not limited to either residual or solution error minimization; one could also combine these two errors into a hybrid loss function, or use quantities of interest like energy instead of the solution itself. For example, List \& Thuerey \cite{list2022a} propose a loss function composed of several terms, and Beck \& Kurz \cite{beck2021} employ the energy spectrum in the loss function. Such hybrid methods have the potential to bridge the common discrepancy between functional and structural approaches (see, e.g.,\ \cite{prakash2022}).

With $\G$ in closure model form \eqref{eqn:closure_model_form}, the residual minimization of the a posteriori formulation \eqref{eqn:Ltf} resembles the a priori formulation \eqref{eqn:Lot_1}. In contrast to a priori learning, the a posteriori approach aims to learn the term $\m_{\pcal}(\qred;\pphys)$ \textit{such that the solution $\qred$ obtained when solving the reduced model \eqref{eqn:Ltf} is close to the reference trajectory $\qbar$.}

A posteriori learning is a relatively recent technique in the closure modeling field, and less literature has appeared in comparison to a priori learning\footnote{note that in PDE-constrained optimization this approach has been in use for much longer \cite{hinze2009}}. 
Early examples of a posteriori training of neural networks in PDEs (mainly fluid dynamics) are the direct field-inversion approach of Duraisamy and collaborators \cite{holland2019} (in the context of RANS), and the Deep learning PDE Model (DPM) of Sirignano et al.\ \cite{sirignano2020} (in the context of LES), which solves the adjoint equations to compute the gradients of the loss function (also known as ``sensitivities''). In \cite{macart2021}, the DPM is extended to train on mean statistics instead of the solution field itself, and efforts are made to interpret the learned neural network closure, as well as its dependence on discretization and the grid. Similarly, \cite{list2022a} learns the sub-grid-scale stress tensor, but uses automatic differentiation to obtain gradients and proposes to use a loss function that contains several terms, reflecting errors in the solution, in the stress tensor, and in the energy spectrum. 
In \cite{shankar2022}, several reduced model forms are considered, ranging from neural ODEs (ResNets) to models that learn a single parameter (the Smagorinsky coefficient); this was extended in \cite{shankar2023} where a neural network model is used as additive correction to the learned Smagorinsky model. 
In \cite{frezat2022}, a posteriori learning is used to learn the stress tensor in quasi-geostrophic turbulence, which is a particularly challenging application because it exhibits a strong energy transfer from small to large scales. 

An ongoing theme in a posteriori learning is how many time steps are ``unrolled'' in evaluating and back-propagating the loss function. Unrolling too few time steps gives only limited gains over a priori learning, while unrolling too many time steps is computationally expensive, has the danger of exploding or vanishing gradients, and can be unrealistic given that turbulent flows are chaotic, meaning that initially close trajectories are expected to diverge \cite{melchers2023,lippe2023,ahmed2023a}. One strategy is to unroll predictions for a small number of the time steps and then only back-propagate to the last time step; this is called the pushforward method \cite{gonzalez2023}. Recently, a similar a posteriori learning approach has appeared in the ROM literature, in which similar challenges are faced \cite{hartmann2021,uy2023}. 

In summary, the main advantage of the a posteriori approach (with solution error minimization) is that one directly targets the accurate approximation of $\qbar$, which is typically the quantity one is interested in. This has been shown to improve stability compared to a priori learning. In addition, if the form of $\G$ does not involve $\F$, like in the evolution form \eqref{eqn:evolution_form}, access to $\F$ is not required, and only reduced solution trajectories $\qbar$ are needed.  The approach implicitly corrects for spatial and/or temporal discretization errors, which can be desirable but can also limit application to different grids or time steps (further discussed in Section \ref{sec:discretization}). The main disadvantage of the a posteriori approach is its computational expense, and depending on the precise form of the reduced model, it requires differentiable solvers (typically using adjoints or automatic differentiation) and a judicious choice of the number of unrolled time steps in the loss function.

\subsection{Approach 3: Reconstruction}\label{sec:reconstruction}


A third approach, which can be used in both a priori or a posteriori fashion, is the reconstruction of the full solution field $\q$ (given the reduced field $\qbar$) and then applying the original model $\F$. So, instead of trying to minimize the error in the \textit{reduced} solution (as given by \eqref{eqn:Ltf}), the goal is to minimize the error in the \textit{reconstructed} solution:
\begin{equation}\label{eqn:closure_problem_rec}
\pcal^{*} = \argmin_{\pcal} \Jr_{\pcal}(\qbar,\q) 
\end{equation}
where the loss function typically involves a solution error 
\begin{equation}\label{eqn:objective_function_rec}
    \Jr_{\pcal}(\qbar,\q) = \| \R_{\pcal}(\qbar) - \q \|^2,
\end{equation}
or a residual error, e.g.,
\begin{equation}\label{eqn:objective_function_rec_res}
    \Jr_{\pcal}(\qbar,\q) = \| \Fbar(\R_{\pcal} (\qbar)) - \Fbar(\q) \| = \| \dd{\qbar}{t} + \bar{\f}(\R_{\pcal}(\qbar)) \|^2,
\end{equation}
and the goal is to find the parameters of the reconstruction operator $\mathcal{R}_{\pcal}$. Once $\pcal$ is determined, one solves
\begin{equation}
    \Fbar(\R_{\pcal}(\qred);\pphys) = 0. 
\end{equation}
We call this approach \textbf{reconstruction}, although the term deconvolution is also commonly used \cite{berselli2006a,stolz1999}. The main challenge is obviously to determine $\R$, which is an inverse problem that is generally ill-posed. In addition, the reconstruction operator is a non-local operator \cite{pawar2020}. Classically, the reconstruction operator has been determined using approaches like the iterative van Cittert method \cite{sagaut2006}. Recently, machine learning alternatives have been proposed, in which $\R$ is modeled by a neural network \cite{maulik2017,maulik2018,yuan2020,yuan2021,fukami2021}. Neural networks are promising because identifying the reconstruction operator is an ill-posed problem that can benefit from data-driven approaches. These methods can be considered an improvement over approximate deconvolution techniques because they do not need the assumption of an invertible filter, but open challenges include how such data-driven methods can generalize across physical regimes, numerical schemes, and discretizations.


With objective function \eqref{eqn:closure_problem_rec}, reconstruction is a special case of a priori learning. By replacing $\qbar$ with $\qred$ in \eqref{eqn:closure_problem_rec} and \eqref{eqn:objective_function_rec}, the approach can also be formulated in terms of a posteriori learning, but this is not common in literature \cite{maulik2017,maulik2018,yuan2020,yuan2021}. Compared to formulation \eqref{eqn:closure_problem_tf}, there is less guarantee that the solution $\qred$ will be close to $\qbar$. Additionally, the use of a posteriori reconstruction operator learning is limiting because of dramatic increases in memory requirements to reconstruct the full flow-field as an indirect route to computing subgrid stresses. We note that the objective functions \eqref{eqn:objective_function_rec} and \eqref{eqn:objective_function_rec_res} are similar to those used when training autoencoders \cite{goodfellow2016}. 

\subsection{From offline to online learning}
The a priori and a posteriori learning approaches from Sections \ref{sec:operator_fitting} and \ref{sec:solution_fitting} are generally used in an offline setting, where the term ``offline'' is used here to indicate that the high-fidelity model is not required in the learning phase. Once training data has been generated by running the high-fidelity model at different parameter values, the high-fidelity model is not needed anymore, and the reduced model can be learned (in a priori or a posteriori fashion). In contrast, recent research is being directed at \textit{online} learning, in which one includes the high-fidelity model in the learning phase. For example, Rasp et al.\ \cite{rasp2020} propose to run the high-fidelity model in parallel with the training process so that it can nudge the reduced model into the right direction and prevent instabilities. Eventually, the idea is that the trained model becomes accurate and stable enough so that it can be run without this coupling. Another example is the dynamic deep learning closure method \cite{sirignano2023}, in which high-fidelity simulations are run in parallel with the reduced model, but only on a small part of the simulation domain. The closure model parameters are dynamically adjusted based on the difference between the high-fidelity model and reduced model on this subdomain, and then used throughout the entire domain. With this approach, no DNS training data needs to be stored, and there is no issue with extrapolation in time or in parameters because the data used for learning the closure model is fit for a purpose. 

We note that the concept of online learning is closely related to the setting in which reduced model parameters (and state) are continuously updated through (observation) data that arrives sequentially in time---we refer to this as data assimilation (see Section \ref{sec:DA}). An example is a neural network-based closure model that is enhanced by using data assimilation techniques \cite{pawar2021}. In addition, another related concept is that of \textit{active} learning \cite{e2020d,fasel2022}, which is typically used to indicate that the samples of the high-fidelity model are iteratively chosen, based for example on the accuracy of the reduced model---see \cite{maddu2023,zhao2018}.

\subsection{Reinforcement learning}\label{sec:RL}
The objective functions in Sections \ref{sec:operator_fitting}-\ref{sec:reconstruction} fit within the framework of supervised learning: the high-fidelity model provides labeled data, and the purpose of the reduced model is to reproduce this accurately. Some limitations of the supervised approach are the need for generating training data (expensive), the issue of defining the reduction operator and reduced training data (the reduction operator is, in practice, often implied by the discretization and therefore not fully known---see Section \ref{sec:discretization}), and the issue of stability, robustness, and the possible need for differentiable programs (see Section \ref{sec:operator_fitting}). 
Reinforcement learning promises to circumvent such issues \cite{novati2021automating,bae2022,kurz2022,kim2022,beck2023toward}: it trains an agent while it is interacting with a dynamical environment. This agent is typically a neural network and the environment is the reduced model, and the agent acts to determine a policy (e.g.,\ the Smagorinsky coefficient, see \ref{sec:LES}) such that a certain reward (objective function) is optimized. For recent review papers, see \cite{garnier2021,viquerat2022}. Some listed advantages of the approach are that only statistics and no high-fidelity data are required, it does not employ adjoints or differentiable solvers (but approximates gradients of the cost function through sampling), and it works directly in the fully discrete reduced model to account for both the closure model error as well as for discretization errors \cite{bae2022}. \RC{rev3_comment28} \revthree{Compared to earlier sections, reinforcement learning shares most similarity with a posteriori learning.} However, one should realize that with a posteriori learning, training based on statistics is also possible, as well as accounting for the discretization error, and gradients can be computed more accurately.




\section{Physics-constrained learning}\label{sec:physics-constrained}

The reduced model forms for $\G$ discussed in Section \ref{sec:problem_formulation} differed in the degree to which ``known physics'' was included. Depending on the application, it is often possible to further inject physics knowledge into the form of $\G$ or in the objective function. In general, the ability to respect physical laws such as conservation and invariance are of key importance for the broad-spectrum utilization of machine learning models for computational physics applications. The development of data-driven models that adhere to such laws leads to improved interpretability and explainability, generalization and robustness, and data efficiency \cite{karniadakis2021physics,meng2022physics}. Therefore, there is a concerted effort to develop algorithms that have constraints embedded in them to avoid violating physical laws. Constraints can be embedded in such data-driven models in a ``soft'' manner where data augmentations or weak penalties in loss functions are used to approximately satisfy constraints. A popular example of such a model is a PINN, mentioned in Section \ref{sec:problem_formulation}, that penalizes deviations from governing laws as a residual term in the model optimization \cite{karniadakis2021physics,cai2021,raissi2019}. These models are typically easier to construct and can be readily applied to existing workflows. In contrast, ``hard'' constraints require the construction of machine learning models that can satisfy constraints by design. While these algorithms require greater turnaround time for construction, training, and deployment, they can guarantee constraint satisfaction during interpolation as well as extrapolation. Examples of model development that constrain a data-driven model by design to permitted behavior include symmetry preservation in graph neural networks (such as for molecular dynamics \cite{e2023} and rotationally invariant flows \cite{shankar2023importance}). \revone{Figure~\ref{fig:soft_hard} shows an example of imposing mass conservation for a 2D incompressible flow as a soft constraint (in a PINN-like fashion), compared to hard constraint (through the definition of streamfunction). Although the hard constraint approach satisfies the divergence-free condition by construction (both during training and deployment), the prediction of the streamfunction is a non-trivial task by itself.} In the following sections, we explore certain examples of soft and hard constraint-based machine learning with applications to closure modeling.


\begingroup
\begin{figure}[hbtp]
\fontfamily{lmss} 
\fontsize{9pt}{9pt}\selectfont
\centering 
\def\svgwidth{\textwidth}
\begingroup%
  \makeatletter%
  \providecommand\color[2][]{%
    \errmessage{(Inkscape) Color is used for the text in Inkscape, but the package 'color.sty' is not loaded}%
    \renewcommand\color[2][]{}%
  }%
  \providecommand\transparent[1]{%
    \errmessage{(Inkscape) Transparency is used (non-zero) for the text in Inkscape, but the package 'transparent.sty' is not loaded}%
    \renewcommand\transparent[1]{}%
  }%
  \providecommand\rotatebox[2]{#2}%
  \newcommand*\fsize{\dimexpr\f@size pt\relax}%
  \newcommand*\lineheight[1]{\fontsize{\fsize}{#1\fsize}\selectfont}%
  \ifx\svgwidth\undefined%
    \setlength{\unitlength}{990bp}%
    \ifx\svgscale\undefined%
      \relax%
    \else%
      \setlength{\unitlength}{\unitlength * \real{\svgscale}}%
    \fi%
  \else%
    \setlength{\unitlength}{\svgwidth}%
  \fi%
  \global\let\svgwidth\undefined%
  \global\let\svgscale\undefined%
  \makeatother%
  \begin{picture}(1,0.26742424)%
    \lineheight{1}%
    \setlength\tabcolsep{0pt}%
    \put(0,0){\includegraphics[width=\unitlength,page=1]{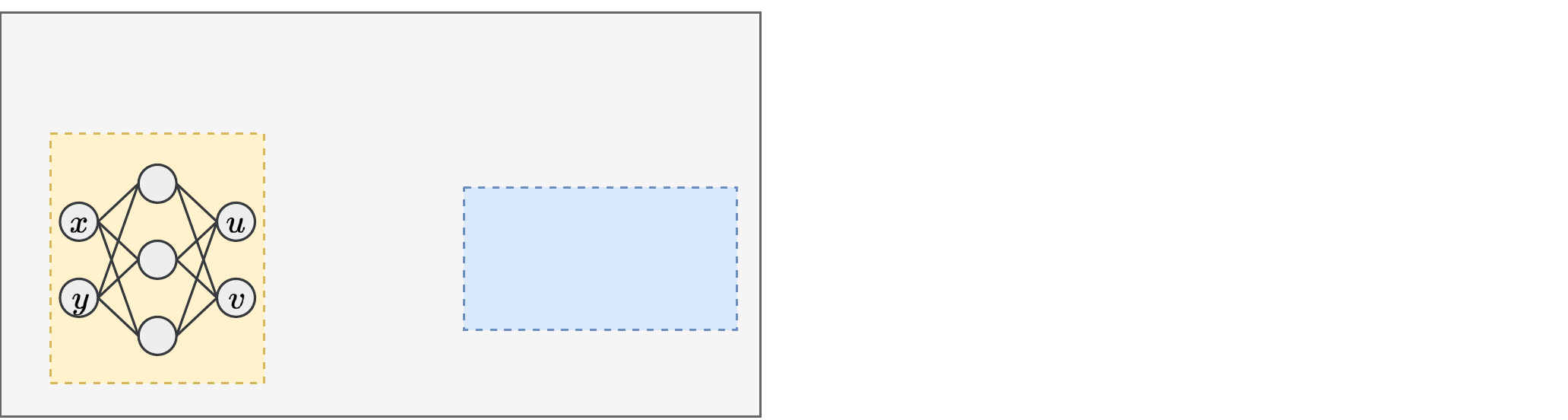}}%
    \put(0.04698323,0.23468868){\makebox(0,0)[lt]{\lineheight{1.25}\smash{\begin{tabular}[t]{l}Mass conservation as a soft constraint\end{tabular}}}}%
    \put(0,0){\includegraphics[width=\unitlength,page=2]{soft_hard_constraint_v2.pdf}}%
    \put(0.84351058,0.09774294){\makebox(0,0)[lt]{\lineheight{1.25}\smash{\begin{tabular}[t]{l}$\dd{u}{x} + \dd{v}{y}=0$\end{tabular}}}}%
    \put(0,0){\includegraphics[width=\unitlength,page=3]{soft_hard_constraint_v2.pdf}}%
    \put(0.30360191,0.09774294){\makebox(0,0)[lt]{\lineheight{1.25}\smash{\begin{tabular}[t]{l}$\mathcal{L}=\|\dd{u}{x} + \dd{v}{y}\|^2$\end{tabular}}}}%
    \put(0.560163,0.23478733){\makebox(0,0)[lt]{\lineheight{1.25}\smash{\begin{tabular}[t]{l}Mass conservation as a hard constraint\end{tabular}}}}%
  \end{picture}%
\endgroup%
 
\caption{\revone{An example of imposing mass conservation for a 2D incompressible flow as a soft constraint (left) and hard constraint (right). For the hard constraint, the streamfunction $\psi$ is learned instead of the velocity field; the velocity field follows as $\vt{u} = \nabla \times \psi$, which by construction satisfies $\nabla \cdot \vt{u}= 0$.
\label{fig:soft_hard}}}
\end{figure}
\endgroup

\subsection{Soft constraints}\label{sec:soft_constraints}

The easy application of soft constraints to machine learning models makes their use attractive for several turbulence closure modeling strategies. \RC{rev3_comment29}\revthree{In the scientific machine learning community, PINNs \cite{raissi2019} are probably the most popular approach for training neural networks using soft constraints, originally introduced by Lagaris et al.\ \cite{lagaris1998artificial}.} These constraints are typically in the form of loss-function regularizations that comprise governing laws. Eivazi et al., \cite{eivazi2022physics}, proposed the use of PINNs for solving the Reynolds-averaged Navier-Stokes equations. In this study, the authors were able to bypass the use of specific model forms for the Reynolds stresses by simply using a governing equation residual loss.  A similar study is seen in \cite{hanrahan2023studying}, where sparse experimental data is also assimilated in the process of training a PINN. Ref.\ \cite{kag2022physics} introduces a two-model formulation where PINN models are trained for low- and high-wavenumber information to obtain more accurate predictions for two-dimensional turbulence. It was observed that a dedicated architecture for high-wavenumber information aided in generalization, \RC{rev3_comment30}\revthree{circumventing the spectral bias problem \cite{rahaman2019spectral}}. In \cite{guan2023}, a hybrid formulation is used where a closure is constructed with enstrophy conservation enforced by a loss-term regularization and where rotational equivariance is preserved by a group-equivariant convolutional neural network. The formulation is applied to construct a closure for a two-dimensional turbulent flow. In \cite{zhao2024mitigating}, a loss-function augmentation is used to restrict the evolution of a machine learning closed system to a low-dimensional phase space to avoid a posteriori errors. This augmentation improved the quality of an a priori trained closure without increased propagating error during a posteriori deployment. \RC{rev2_comment6}\revtwo{Finally, another approach to satisfying a class of constraints weakly is through the use of data augmentation. Here, augmentations to training datasets can be performed to implicitly inform the learning algorithm about invariances or equivariances in the function approximation. For example in \cite{brandstetter2022lie}, data augmentation with Lie symmetries not only leads to the improved satisfaction of the same but also improves sample complexity. Similar operations were also deployed in \cite{kim2020deep} for generating synthetic turbulence snapshots. In \cite{subel2021data}, a similar data augmentation was performed for learning subgrid stress tensors in two-dimensional turbulence with improved results compared to the standard dataset.}

\subsection{Hard constraints}\label{sec:hard_constraints}

Hard constraints are typically application specific, and we focus on turbulence closure modeling. It is well known that subgrid stress tensors in LES, or Reynolds stress tensors in RANS, obey certain laws such as aggregate dissipative behavior, symmetry, invariance, etc.\ \cite{speziale1985,silvis2017}. Therefore, it is desirable to construct data-driven prediction techniques for these quantities that adhere to these constraints during inference. The simplest efforts to constrain predictions from data-driven closures have involved ad-hoc truncations that preserve numerical stability in deployments. These typically occur via clipping to ensure guaranteed dissipation of the subgrid stress \cite{maulik2019,beck2019}. Other efforts include constrained predictions via projection of data-driven predictions onto the so-called realizability triangle \cite{wang2017comprehensive}. Another popular approach for constructing data-driven approximations of such tensors are via tensor-basis neural networks (TBNNs) \cite{ling2016,ling2016a,wang2017comprehensive,bose2024}, which guarantee a Reynolds stress prediction as a linear expansion of strain and rotation tensors following Pope's generalized eddy viscosity hypothesis \cite{pope1975more}. \RC{rev1_comment1}\revone{The stress tensor $\vt{\tau}$ - denoted by $\mdiv_{\pcal} (\qbar;\pphys)$ in \eqref{eqn:commutator_tensor} - is expressed as
\begin{equation}\label{eqn:TBNN}
    \vt{\tau} = \sum_{i=0}^{10} \alpha_{i} \vt{T}_{i},
\end{equation}
where the coefficients $\alpha_{i}$ are functions of a set of scalar invariants $\lambda_{1}, \ldots, \lambda_{5}$, which in turn depend on the strain-rate tensor $\vt{S}$ and the rotation tensor $\vt{R}:=\frac{1}{2}(\nabla \vt{u} - (\nabla \vt{u})^{T})$, and $\vt{T}$ is the tensor integrity basis, also a function of $\vt{S}$ and $\vt{R}$. For example, $\lambda_{1} = \mathrm{tr} (\vt{S}^2)$ and $\vt{T}_{1} = \vt{S}$. The functions $\alpha_{i} = f_{i}(\lambda_{1}, \ldots, \lambda_{5};\pcal)$ are modelled as neural networks with parameters $\pcal$. In the large eddy simulation context, similar symmetries may also be leveraged for constrained predictions from neural network models \cite{silvis2017,razafindralandy2006,speziale1985,speziale1989}. The important message is, that by choosing the form \eqref{eqn:TBNN}, one has directly encoded the Galilean, rotational, and reflectional invariances of the stress tensor, independent of the neural network architecture or its parameters.}

In several applications, property preservation can be posed as linear constraints. In such cases, neural network architectures that guarantee satisfaction of these constraints have also been proposed \cite{beucler2019,kashinath2020enforcing}. These include the solution of an additional system of equations to project the output of the neural network to a subspace that satisfies these equations. However, it must be noted that this involves the inversion of a matrix that grows with the size of the output for each forward pass of the neural network. Another approach pioneered in \cite{vangastelen2023} introduces an additional set of latent variables and specifies an explicit form of the closure (a combination of skew-symmetric and dissipative terms) to guarantee energy conservation and stability of the closure; see also \cite{charalampopoulos2022} for related work. In reduced-order models, similar ideas have been used \cite{mohebujjaman2019,ahmed2021closures,callaham2022a}. 

Finally, a popular approach to building closure models in data-limited scenarios is through field inversion and machine learning \cite{parish2016paradigm}, which allows for the PDE-constrained optimization of defects in turbulence models using experimental or sparse sampling data. If the defects to turbulence models are appropriately defined, constraints can be satisfied by design both during the optimization as well as the deployment procedure \cite{shankar2023differentiable}. A similar procedure may also be followed via reinforcement learning for closure modeling \cite{beck2023toward,novati2021automating,kurz2023deep} (see also Section \ref{sec:RL}), \RC{rev3_comment31} \revthree{which offers the advantage of sparse data requirements,} although computational costs associated with learning are significant in the absence of gradient information. We would like to emphasize that hard and soft constraints are not mutually exclusive and that data-driven modeling approaches can blend both concepts for balancing accuracy, constraint satisfaction, and ease of learning and deployment. For example, in \cite{fiore2022physics}, hard constraints are used for tensor symmetries using TBNNs and soft constraints are imposed for the smoothness of heat fluxes.

In the more general context of symmetry preservation for fluid flow forecasting, several novel neural network algorithms have been proposed that can be repurposed for closure modeling. Wang et al.\ \cite{wang2021} incorporate translational, scaling, and rotational equivariance of the Navier-Stokes equations in the neural network design. In \cite{siddani2021, pawar2023}, novel convolutional neural network architectures that preserve frame invariance are used for fluid-flow predictions. Specifically, Pawar et al.\ \cite{pawar2023} use a novel architecture for improved subgrid stress modeling of two-dimensional large eddy simulation. Guan et al.\ \cite{guan2023} used group convolutional neural networks and demonstrated that the preservation of symmetries alleviated the need for much larger datasets.

Another approach to preserving desirable properties during regression tasks is via the use of non-parametric methods like random forests. In \cite{yuval2020}, a random-forest regressor is used for predicting subgrid stresses for climate models, with guaranteed satisfaction of symmetry preservation. This is because random-forest predictions are obtained by averaging splits of data that come from the convex hull of the training data set. As long as the generated training data sets guarantee a particular symmetry or structure, the prediction of the dependent variable will also satisfy this property by design. Note that for complicated properties that emerge from complex interactions between machine learning predictions and the numerical solver, such an approach may be incomplete.

\section{Discretization aspects}\label{sec:discretization}
The discussion in Sections \ref{sec:problem_formulation} - \ref{sec:physics-constrained} was mostly agnostic of the spatial and temporal discretization methods needed to solve the true model \eqref{eqn:PDE} and approximate model \eqref{eqn:filteredPDE}. Spatial and temporal discretizations give rise to additional levels of approximation, and we indicate by $\Fh$ the (fully discrete) \textit{high-fidelity model}, which satisfies
\begin{equation}\label{eqn:PDE_discrete}
\Fh(\qh^{n};\pphys) = 0,
\end{equation}
where $\qh^{n} \approx \q(\vt{x},t^{n})$. Similarly, $\Gh$ denotes the (fully discrete) \textit{low-fidelity model}, which satisfies
\begin{equation}\label{eqn:filteredPDE_discrete}
    \Gh(\qhred^{n};\pphys) = 0,
\end{equation}
and $\qhred^{n} \approx \qred(\vt{x},t^{n})$. 

\subsection{Spatial discretization and neural ODEs}\label{sec:neuralODE}
In practice, the low-fidelity model is often constructed based on a high-fidelity model that is only spatially discretized, leaving time continuous. An important example is when the high-fidelity model is in the form of equation \eqref{eqn:evolution_form}, leading to 
\begin{equation}\label{eqn:discrete_closure_model_form_ddt}
    \Gh(\qhred;\pphys) = \frac{\rd \qhred}{\rd t} + \fh(\qhred;\pphys) + \mh(\qhred;\pphys).
\end{equation}
In this case, the spatial discretization $\fh$ of the high-fidelity model can be reused but on a coarser grid, and existing high-fidelity codes can simply be extended by adding a closure term. Instead of a coarser grid, a low-order method can also be used, like in \cite{kang2023}. Furthermore, this notation also covers the case where $\qhred$ is being represented by the coefficients of an ROM, see, e.g.,\ \cite{gupta2021a,ahmed2021closures}. 

In scientific machine learning, the common choice is to represent $\mh$ by a neural network, resulting in so-called neural (closure) ODE results \cite{gupta2021a,melchers2023,chen2019,kidger2022} (the concept has been denoted as the Universal Differential Equation in \cite{rackauckas2021}). Casting the closure problem as a neural ODE has the benefit that one can formulate the neural network training as an optimal control problem, which can be solved for example by using adjoint methods, as proposed in \cite{chen2019}, and the parameters of the network can be interpreted as the control variables. 
One difference is that in the classical optimal control setting \cite{pontryagin1986,kirk2004} the control variables are time-dependent, whereas the neural network parameters in \cite{chen2019} are time-independent. These are important differences, some of which are addressed in \cite{massaroli2020}. A recent approach in which the closure model parameters change dynamically over time is given in \cite{sirignano2023}.

Using neural networks for closure modeling of form \eqref{eqn:discrete_closure_model_form_ddt} is advantageous because time is continuous, so the neural networks are expected to also work when different time integrators or different time steps are used. However, the fact that one is dealing with a spatially discretized system means that the neural network is linked to the particular discretization method and grid that were present in the training data, and will not easily generalize to different grids or different discretization methods. On the other hand, an advantage of learning discretized operators is that they can correct for discretization errors \cite{sirignano2020}. 

As mentioned, the time-continuous nature of neural ODEs allows one to formulate the adjoint equations associated with the optimization problem \eqref{eqn:Ltf}, see, e.g.,\ \cite{chen2019}, which can be solved backward in time to yield the sensitivities ($\dd{\mathcal{L}}{\theta}$). Alternatively, one could formulate the time-discrete equations and then determine sensitivities by backpropagating through the entire ODE solver. Optimize-then-discretize or discretize-then-optimize is a topic that has been discussed in the optimization literature---see textbooks like \cite{hinze2009, gunzburger2003}. \revone{Figure~\ref{fig:optimize_discretize} shows a pictorial example of both paradigms.} Recent research suggests that there are some important benefits to \RC{rev2_comment7}\revtwo{discretize-then-optimize} \cite{onken2020,melchers2023}, like obtaining the exact gradient (of the time-approximate solution) instead of the approximate gradient that is returned by adjoint methods.



\begingroup
\begin{figure}[hbtp]
\fontfamily{lmss} 
\fontsize{9pt}{9pt}\selectfont
\centering 
\def\svgwidth{\textwidth}
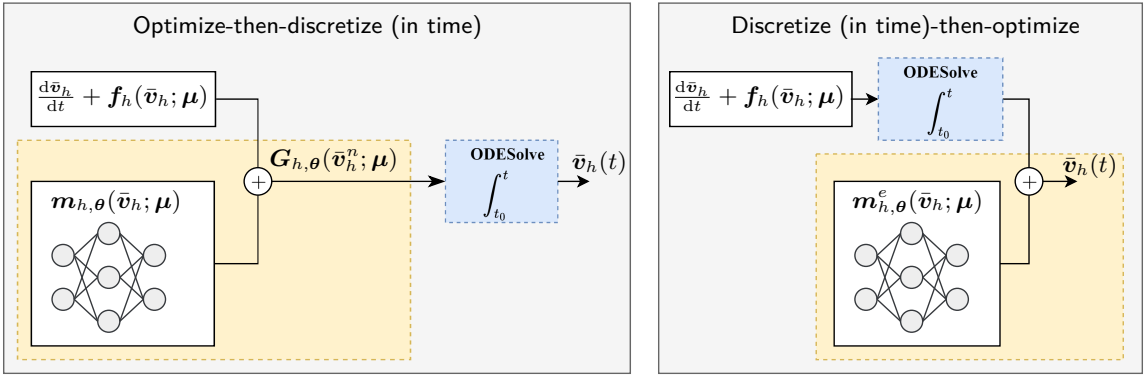 
\par
\caption{\revone{Optimize-then-discretize vs.\ discretize-then-optimize approaches. In orange: the optimization of a neural network to learn a time-continuous closure (left) or a time-discrete correction (right). In blue: the time integration/discretization step.
\label{fig:optimize_discretize}}
}
\end{figure}
\endgroup

\subsection{Temporal discretization and autoregressive methods}
An example of a ``\revtwo{discretize-then-optimize}'' formulation is 
\begin{equation}\label{eqn:discretized_evolution}
    \qhred^{n+1} - \qhred^{n} + \Delta t \left( \fh(\qhred^{n};\pphys) + \meh(\qhred^{n};\pphys)\right) = 0,
\end{equation}
which corresponds to \eqref{eqn:filteredPDE_discrete} when the discretized evolution form is taken for $\Gh$. Although there are several advantages of discrete methods (exact gradients, correction for discretization errors), one important disadvantage is that the closure model is not only linked to the spatial discretization but also to the time step used during training. 

Beyond the closure forms \eqref{eqn:closure_model_form} and \eqref{eqn:closure_model_form_ddt}, fully discrete approaches are quite commonly used in reduced models of the evolution form \eqref{eqn:neuralODE}, in particular so-called autoregressive methods. Although these fall strictly speaking beyond the concept of closure modeling, they deserve to be mentioned given their recent success in weather prediction \cite{pathak2022,lam2023} and in general in their use in predicting chaotic dynamical systems \cite{vlachas2020,vlachas2022}. The main idea is to generalize the discrete version of \eqref{eqn:neuralODE}, e.g.,\ $\qhred^{n+1} = \qhred^{n} - \Delta t \me(\qhred^{n};\pphys)$, to the form 
\begin{equation}\label{eqn:autoregressive}
    \qhred^{n+1} = \mathcal{N}_{\pcal}(\qhred^{n},\qhred^{n-1},\qhred^{n-2},\ldots;\pphys).
\end{equation}
In other words, one learns an operator $\mathcal{N}_{\pcal}$ that uses past solution states to perform time stepping. Even though the closure form \eqref{eqn:discretized_evolution} intuitively seems easier to learn because it involves only a ``correction'' to the known physics, its mathematical form is less expressive, and it has been observed that this limits the model performance \cite{kochkov2021}. In fact, \eqref{eqn:autoregressive} has a non-Markovian form, allowing it to incorporate history, which has its foundations in the Mori-Zwanzig theory, as will be discussed in Section \ref{sec:MZ}. 
The loss function typically involves the error in the solution, like in Section \ref{sec:solution_fitting}, so this is a type of a posteriori learning. Given that multiple past states are included, the computation of the derivatives of the loss function through back-propagation will be more involved than for a Markovian approach like \eqref{eqn:discretized_evolution}. In addition, one would expect again that the learned operator is only accurate for particular time steps. However, the impressive performance of weather predictions \cite{pathak2022,lam2023,bodnar2024} indicates that this is perhaps less of a concern for practical problems. 
An example of an autoregressive method for LES is given in \cite{li2022a}, where an FNO is employed. Autoregressive methods have also been used in message-passing PDE solvers \cite{brandstetter2023}. 

\RC{rev3_comment33} \revthree{Another example where form \eqref{eqn:autoregressive} is employed, is in `model-free' machine learning for forecasting of chaotic systems \cite{pathak2018a,vlachas2020,vlachas2022}. In that case, $\qhred$ typically consists of an internal hidden memory state and an observable, which are evolved in time with a recurrent neural network \cite{vlachas2022}. The foundation behind this approach lies in Takens' theorem \cite{takens1981}, which states that a dynamical system can be identified from partial observations of the state, as long as sufficient history (delays) of these partial observations is available. This is reminiscent of the Mori-Zwanzig formalism, which offers a framework to reduce the dimension a dynamical system by including history effects. However, Takens' theorem instead aims at reconstructing the full state based on time-delayed observations.}



\subsection{Learned discretizations}
An approach that is related to the two previous sections entails directly learning a (coarse) discretization scheme that implicitly includes the effect of the closure term, such as the ``learned interpolation'' method \cite{bar-sinai2019}. So, instead of employing the spatial discretization of the high-fidelity model on a coarse mesh and adding a correction term to it, cf.\ \eqref{eqn:discrete_closure_model_form_ddt}, one learns a discretization scheme (an improved interpolation method based on neural networks) for a coarse grid that effectively takes into account unresolved physics. This approach is rather intrusive because it requires one to change an existing code at the level of the discretization, but comes with the advantage that the learned interpolation (effective discretization stencil) is in some sense interpretable. Interestingly, the implicit LES formulation \cite{grinstein2007} (which was developed long before machine learning came into play) also falls within this category---it does not feature a closure model, but relies on the fact that the discretization typically introduces additional dissipation compared to the true physics, which mimics the additional dissipation that is usually provided by the closure model.

A less intrusive approach is the ``learned correction'' approach \cite{um2020,kochkov2021,dresdner2023}, which is simpler to implement and more flexible. After performing a time step with the coarse grid operator, the solution is corrected with a learned term. This approach is in fact similar to the fully discrete formulation \eqref{eqn:discretized_evolution} when the term $\Delta t \meh$ is identified as a ``learned correction''. In \cite{stachenfeld2022a}, a similar approach is taken, but the known physics part expressed by $\fh$ is not included, so the neural network has to learn the entire effect of both known physics and correction.

It should be noted that in many practical LES models, the filtering operation \eqref{eqn:continuous_filter} is not explicitly defined, but implied by the discretization \cite{beck2023toward}. The alternative, namely explicit filtering, has the important advantage that it converges to the ``true filtered'' equations upon grid refinement \cite{Lund2003,Gullbrand2003}. However, because this is generally considered to be computationally too expensive, the implied filter approach is often preferred. The disadvantage, however, is that it is difficult to define filtered reference data (needed to train neural networks) because the filter is not explicitly known. 


\subsection{Discretization invariance}
All approaches that are spatially or temporally discrete suffer to some extent from the issue that generalization to different grids, time steps, and discretization methods is difficult. Neural operators like DeepONet and FNO are an important step in the direction of methods that are grid-independent because they are mappings between (infinite-dimensional) function spaces \cite{li2021,wang2022,lu2021,kovachki2023}. In practice, training data for such operators still involves discrete input and output data. Recently, a class of neural operators has been proposed in which the discrete input-output representation and their continuous function space realizations are equivalent \cite{bartolucci2023}. 



\section{Connections to inverse problems, PDE-constrained optimization, and statistical methods}\label{sec:inverse_etc}

\subsection{Inverse problem specification}\label{sec:inverse}
The minimization problems \eqref{eqn:closure_problem_of}, \eqref{eqn:closure_problem_tf} and \eqref{eqn:closure_problem_rec} from Section \ref{sec:objective_functions} can be interpreted as inverse problems. They are similar in that some output quantity is known based on the high-fidelity simulation (e.g.,\ the exact reduced solution or the exact reduced operator), and one can try to find a model that best approximates that output. The inverse-problem nature implies that a unique solution typically does not exist and that solutions can be sensitive to perturbations.

For example, the a posteriori learning problem can be formulated in terms of an inverse problem, which typically features two components: a forward model, 
\begin{equation}\label{eqn:forward_model}
  \G_{\pcal}(\qred;\pphys) = 0,
\end{equation}
and an observation equation,
\begin{equation}\label{eqn:observation_equation}
    \obs = \vt{h}_{\pcal}(\qred) + \varepsilon,
\end{equation}
where $\obs$ denotes the observations and $\vt{h}$ denotes the observation operator. In inverse problem terms, the art is to find the parameters $\pcal$ given the observations and the forward model. In the closure model framework outlined in Sections \ref{sec:problem_formulation} and \ref{sec:objective_functions}, the observations come from high-fidelity simulations and the observation operator extracts the corresponding data from the reduced model. For example:
\begin{align}
  &\textrm{observations of the state: }   & \obs &= \qbar, \qquad \vt{h}_{\pcal}(\qred)  = \qred,  \\
  &\textrm{observations of the closure term: }   & \obs &= \C(\q,\qbar;\pphys), \qquad \vt{h}_{\pcal}(\qred)  = \m_{\pcal} (\qred;\pphys).
\end{align}
The a priori learning problem has a similar form, with the important difference that the observation equation \eqref{eqn:observation_equation} changes to $\obs = \vt{h}_{\pcal}(\qbar)  + \varepsilon$, which can be solved \textit{independently} from the forward model \eqref{eqn:forward_model} (see Figure \ref{fig:a-priori_a-posteriori}).

The formulation above can aid in employing concepts from the inverse problem community, like regularization, that address the ill-posedness of the learning problem \cite{hansen2010,stuart2010,calvetti2018}. For example, a regularized inverse problem could read
\begin{equation}\label{eqn:inverse_formulation}
    \min_{\theta} \| \obs -\vt{h}_{\pcal}(\qred) \|^2 + \alpha T(\pcal),
\end{equation}
with a regularization function $T$, e.g.,\ $T = \| \pcal \|^2$, and a hyperparameter $\alpha$. In addition, the connection between machine learning in general and inverse problems has been the focus of many efforts. For instance, derivative-free optimization of neural network parameters (i.e., without back-propagation) has been achieved by means of a modified ensemble Kalman filter \cite{haber2018never} and ensemble Kalman inversion \cite{kovachki2019ensemble}. This can potentially mitigate the need for differentiable solvers in a posterior learning. Data-driven parameterization studies (in the context of weather and climate predictions) have benefited from the ensemble Kalman inversion (EKI) approach. Lopez-Gomez et al.\ \cite{lopez2022training} posed the learning of SGS turbulence and convection as a regularized inverse problem and used EKI and unscented Kalman inversion to solve it. In \cite{pahlavan2024explainable}, EKI is used to re-train two layers of a previously offline-trained neural network model for gravity wave forcing.



\subsection{PDE-constrained optimization}\label{sec:optimalcontrol}
Because the minimization formulation of the inverse problem, equation \eqref{eqn:inverse_formulation}, is constrained to the (generally) nonlinear forward model dynamics---which is typically a PDE---it is a PDE-constrained optimization problem \cite{hinze2009,gunzburger2003}. This viewpoint, and the realization that the dimension of the parameter space is typically very large, naturally leads to the consideration of adjoint methods that compute the gradient $\dd{J_{\pcal}}{\pcal}$ of the objective function $J_{\pcal}=\| \obs -\vt{h}_{\pcal}(\qred) \|^2$. Formulating the adjoint equations, discretizing them, and solving the optimization problem is the DPM approach taken in \cite{sirignano2023}, for example. Alternatively, one can first discretize the PDE and use back-propagation (reverse-mode automatic differentiation) through the entire PDE solver to get gradients of the cost function with respect to the parameters \cite{list2022a,melchers2023} (see Section \ref{sec:solution_fitting}).

\subsection{Data assimilation}\label{sec:DA}
Another useful perspective is to view the closure problem from the angle of data assimilation \cite{asch2016,ahmed2021closures}. In variational data assimilation, the formulation is similar to equations \eqref{eqn:forward_model} and \eqref{eqn:observation_equation} (but formulated for the discrete (finite-dimensional) system), with several main differences: observations are not high-fidelity simulations but actual measurements; the parameters are not neural network parameters but initial conditions; and the optimization problem is split into sequential steps, corresponding to temporal windows for which the measurement data are becoming available. 

The data assimilation perspective is useful because it provides a natural approach to extend the learning of closure models by high-fidelity data to include (streaming) measurement data. The parameters of closure models can be seen as time-dependent quantities that are changing while high-fidelity or measurement data arrives, in line with the online learning approaches of \cite{rasp2020} and \cite{pawar2021}, \RC{rev3_comment35} and the reinforcement learning approach explained in \ref{sec:RL}. The data assimilation perspective and PDE-constrained optimization formulation naturally connect to the concept of optimal control. In optimal control, one seeks (time-dependent) parameters that minimize an objective function, subject to given dynamics.

\subsection{Statistical approaches}\label{sec:statistical}
The approaches in Sections \ref{sec:inverse}, \ref{sec:optimalcontrol}, and \ref{sec:DA} are variational in the sense that a minimization problem is solved. They provide a single-point estimate of the closure term. This is the most common approach in the closure modeling community. A logical alternative is to consider statistical approaches like Bayesian methods, which can naturally include prior knowledge of the parameters. Bayesian methods treat the parameters of interest as random variables and provide a full probabilistic description of the knowledge about the learned parameters. Such methods have the advantage that they naturally come with an expression for the uncertainty in the parameters. However, their computational expense (typically due to sampling the posterior) is high. In the data assimilation community, the use of statistical methods instead of variational methods is in fact also common \cite{asch2016}. Before the rise of scientific machine learning methods, Bayesian methods have been applied to learn parameters of closure models, in particular of RANS model coefficients \cite{edeling2014,xiao2016}. 

In the context of scientific machine learning, the parameters of interest in the Bayesian framework can be either the neural network weights $\pcal$, as in Bayesian neural networks \cite{goan2020bayesian,jospin2022hands}, or the model output itself (e.g., the closure term $\m_{\pcal}$). In \cite{bhouri2022history}, a Bayesian inference framework is developed to learn the distribution of neural network weights, with application to SGS parameterization in the Lorenz 96 model. The statistics of the predicted parameterizations can be approximated by sampling from the posterior distribution of the network parameters. A simpler way to estimate the uncertainty in the model's prediction is through deep ensemble methods \cite{lakshminarayanan2017simple}, where an ensemble is formed by training the same model multiple times with different initializations. Such an approach was employed in \cite{pawar2023,ahmed2023physics} to model the uncertainty in the learned closure model in LES and ROM contexts. Alternatively, the ensemble can be created by varying the model choices (e.g., neural network architectures), training data, the optimization algorithm, dropouts, etc. In addition, machine learning models can be trained to directly predict statistical measures of the closure terms (e.g., mean and variance). Then, the loss function may be written in terms of the log-likelihood of the corresponding distribution. The combination of data assimilation and machine learning in a Bayesian setting has been explored in the more general context of learning dynamics \cite{bocquet2020}. In the closure modeling context, this has only being explored recently \cite{maruyama2021}.


\section{Non-locality and the Mori-Zwanzig viewpoint}\label{sec:MZ}

In general, because closure assumes that certain degrees of freedom will not be represented explicitly (either due to necessity or efficiency), its representation can lead to \textit{non-local} effects for the large-scale or resolved degrees of freedom. As an example, in the case of closure for a time-dependent problem due to a coarse spatial resolution, one may need to use  neighboring (or even further apart) spatio-temporal information about a resolved point when constructing its closure. We review temporal and spatial non-local closures separately, while keeping in mind that they can also appear simultaneously.

\subsection{Non-locality in time}

The Mori--Zwanzig (MZ) formalism is a general framework for treating temporal non-locality. It originated in non-equilibrium statistical mechanics as a way to construct transport equations \cite{zwanzig2001} and was later reformulated as a model reduction (closure) approach \cite{CHK00,CHK3,chorin2007} for systems of ordinary differential equations. It is based on the idea that instead of constructing a reduced model directly at the level of the ordinary differential equations (which are usually nonlinear and thus not easily amenable to projection approaches), one can construct a {\it linear} PDE for each resolved degree of freedom. 

For the sake of concreteness, let us assume like in Section \ref{sec:neuralODE} that we have performed only spatial discretization for a time-dependent PDE and that we have obtained a system of ODEs 
\begin{equation}\label{eqn:spatial_discrete}
\frac{\rd \qh}{\rd t} =  \rh(\qh;\pphys),
\end{equation}
where $\qh \approx \q(\vt{x},t)$ is the vector denoting the spatially discretized solution using resolution $h$, resulting in $M$ degrees of freedom. Also, let $\qh^0 \approx \q(\vt{x},0)$ be the vector denoting the discretization of the initial condition. To avoid confusion for the reader, we note that due to the usual convention of presentation of the MZ formalism, the dynamics of the high-fidelity model are written as in \eqref{eqn:spatial_discrete} instead of the residual form $\frac{\rd \qh}{\rd t} +  \fh(\qh;\pphys) =0$, which we have used so far in the text. Thus, $\rh(\qh;\pphys) = - \fh(\qh;\pphys).$

According to the MZ formalism, the $M$-dimensional system of ODEs \eqref{eqn:spatial_discrete} can be transformed into a system of $M$ linear PDEs, one for each spatial-discretization point \cite{zwanzig2001}. Using semigroup notation, the linear PDEs are given by 
\begin{equation}\label{pdes_linear}
\pd{}{t} e^{tL} \q_{hk}^0=L e^{tL} \q_{hk}^0 \; \; \quad \text{for} \quad \;\; k=1,\ldots,M,
\end{equation}
where
\begin{equation}\label{liouville}
L=\sum_{j \in M } \ro_{hj}(\qh^0) \frac{\partial}{\partial \q_{hj}^0}
\end{equation}
is called the Liouville operator and $e^{tL} \q_{hk}^0 = \q_{hk}(t) $. 

Suppose that the vector of initial conditions can be divided as $\qh^0=(\qbar_h^0,\tilde{\q}_h^0),$ where 
$\qbar_h^0$ is the vector of the resolved (typically large-scale) variables and $\tilde{\q}_h^0$ is the vector of the unresolved (typically small-scale) variables. Let $P$ be an orthogonal projection on the space of functions of $\qbar_h^0$ and $Q=I-P$. 

Equation \eqref{pdes_linear} 
can be rewritten as 
\begin{equation}
\label{mz}
\frac{\partial}{\partial{t}} e^{tL}\q_{hk}^0=
e^{tL}PL\q_{hk}^0+e^{tQL}QL\q_{hk}^0+
\int_0^t e^{(t-s)L}PLe^{sQL}QL\q_{hk}^0 \rd s,
\end{equation}
where we have used Dyson's formula
\begin{equation}
\label{dyson1}
e^{tL}=e^{tQL}+\int_0^t e^{(t-s)L}PLe^{sQL} \rd s.
\end{equation}
Equation \eqref{mz} is the MZ identity. It decomposes the right-hand-side of the PDE system \eqref{pdes_linear} into three terms: (1) a term involving only the values of the resolved degrees of freedom at the current time (usually called the ``Markovian'' term), (2) a term that evolves in the space orthogonal to the resolved degrees of freedom (usually called the ``noise'' term), and (3) a term that accounts for the history of the interaction between resolved and unresolved degrees of freedom (usually called the ``memory'' term). The MZ formalism is based on an exact reformulation of the original system and though almost never used in its full generality, it is a good starting point for approximations based on mathematical, physical, and numerical considerations. Specifically, it is able to address cases with very short, moderate, or very long memories. A detailed presentation of the MZ formalism and its relation to other formalisms developed for specific memory types can be found in \cite{givon2004}. Also, we must note that while we have presented the MZ formalism for the case of a system of ODEs, the formalism has also been developed for systems of stochastic differential equations---see, e.g., \cite{ma2019coarse,hudson2020coarse,zhu2021hypoellipticity}. 


The MZ formalism, through various approximations, has been used to obtain reduced order models and closures for a host of systems (see, e.g.,\ \cite{bernstein2007optimal,stinis2012numerical,li2015memory,lei2016,parish2017,parish2017a,jung2017iterative, pan2018,leimkuhler2022efficient}). However, except for special cases, it is difficult to guarantee the stability of reduced order models. A renormalized version of the MZ formalism has been introduced (inspired by renormalization methods in physics), which has allowed the stabilization of such models \cite{stinis2009,stinis2015,price2019renormalized,price2021optimal}.

\RC{rev2_comment5a}\revtwo{Finally, in recent years, data-driven approaches to MZ have also emerged. These approaches use neural networks to develop frameworks that are inspired or can be mapped to MZ \cite{berner2017,brennan2018, wiewel2019latent,wang2020c, lin2021data,levine2022framework,wang2020recurrent,e2023, dietrich2023learning}. Also, there exist approaches that represent one or more parts of the MZ formalism (i.e.,\ Markovian term, memory, and/or noise). These range from the use of effective machine-learned potentials \cite{wang2019machine,durumeric2023machine}, to neural network representation of the memory \cite{ma2018model, qi2022machine, russo2022machine, bhouri2023a, she2023data, gupta2021a, crommelin2021}, to statistical-informed neural networks in order to simulate the noise \cite{zhu2023learning}.} The success of the approaches highlights the need for significant amounts of training data, the right choice of neural network architecture for the calculation of the various MZ terms, as well as the expected issue of guaranteeing stability for the effective equations for the resolved degrees of freedom.





\subsection{Non-locality in space}

Non-locality in space has been used to model phenomena that are not amenable to modeling, e.g., with integer-order derivatives. Non-locality in space aims to accommodate two types of issues: (1) rough solutions that may not be differentiable, and (2) long-range spatial correlations at the scale where we represent a phenomenon. There are two main approaches to account for non-locality in space---(1) the use of integral terms (of convolutional type), e.g., like in peridynamics, and (2) the use of fractional order derivatives (see e.g. \cite{javili2019peridynamics,kharazmi2020fractional, d2020numerical,suzuki2023fractional, du2020nonlocal}. 

\RC{rev2_comment5b}\revtwo{The use of closure exhibiting non-local spatial features is a topic of active research and includes traditional (see e.g., \cite{held2004nonlocal, zhang2020virtual, liu2023systematic} and machine-learning approaches. The machine-learning approaches range from techniques to estimate the various exponents appearing in the relations representing the non-local terms \cite{song2021variable, seyedi2022data}, to using neural networks to represent the kernels appearing in the non-local expressions \cite{joglekar2023machine, de2023machine, you2021data}, to the use of neural networks to represent the whole non-local term \cite{zhou2021learning, charalampopoulos2022, liu2022investigation, guan2022b, tabe2023priori, han2023equivariant}.} The approaches again highlight the need for accurate training data, the right choice of neural network architecture, and the careful coupling with the rest of the dynamics present in the model.  



\section{Multi-fidelity and multiscale viewpoints}\label{sec:multifidelity}
\subsection{Multi-fidelity}



The process of closure modeling and surrogate modeling almost always involves collecting curated datasets of input-output responses, either through high-fidelity simulations or through experiments (one exception is PINN \cite{raissi2019}). However, in practice, the incurred costs of data collection can restrict the quantity and quality of simulations or experiments. Building a model solely from a small high-fidelity dataset or a large low-fidelity dataset is not generally adequate. Multi-fidelity modeling approaches aim to close this gap by leveraging both low- and high-fidelity components to achieve levels of accuracy and efficiency that are not attainable otherwise. The term ``component'' here refers to the fact that different fidelities might appear in terms of collected data, model assumptions, physical complexity, numerical schemes, etc. For example, Howard et al.\ \cite{howard2023} showed that a multi-fidelity DeepONet can be trained by combining low-fidelity data with knowledge of physics (in terms of the governing PDE) as the high-fidelity portion.

Although the current review focuses on approaches that improve the estimation of the \emph{reduced} quantity $\qbar$ (through closure modeling), we note that the vast majority of literature on multi-fidelity modeling is instead focused on directly estimating the quantity of interest $\q$. In other words, the low-fidelity model provides $\q_{low}$ and the high-fidelity model gives $\q_{high}$, representing the low- and high-fidelity approximations of the same quantity $\q$. Then, frameworks are built such that the majority of computations are delegated to low-fidelity models whose outputs can be corrected using a small number of high-fidelity runs. Multi-fidelity modeling approaches have been particularly beneficial for outer-loop applications (e.g., uncertainty quantification and design optimization)---for a survey of such methods, see \cite{peherstorfer2018survey}. Different studies have considered different variants for the low-fidelity model and how it is embedded in the process of learning high-fidelity outputs. For example, Conti et al.\ \cite{conti2023multi} developed a long short-term memory (LSTM) framework to combine abundant low-fidelity data and limited high-fidelity data, where the low-fidelity data are generated from either a pretrained ML model or using coarser spatial/temporal discretizations. The quantity of interest in the multi-fidelity setting was, however, a low-dimensional vector (e.g., the drag and lift coefficients), which was also the theme in the majority of multi-fidelity surrogate modeling literature \cite{zhang2021multi,lamberti2021multi}. Recently, convolutional neural networks (CNNs) and graph neural networks (GNNs) have alleviated this limitation and allowed for dealing with flow field inputs/outputs efficiently \cite{liao2021multi,mondal2022multi,li2023multi}.

In \cite{conti2024multi}, the full field prediction problem is addressed by projecting the low-fidelity solution onto a common subspace spanned by the leading proper orthogonal decomposition (POD) modes in the form of $\q_{low}(\vt{x},t) = \sum_{i=1}^{N} a_{low,i}(t) \vt{\phi}_{i}(\vt{x})$. The low-fidelity POD coefficients are then passed through LSTM to learn their high-fidelity counterparts, giving a high-fidelity approximation $\q_{high}(\vt{x},t) = \sum_{i=1}^{N} a_{high,i}(t) \vt{\phi}_{i}(\vt{x})$. Because POD-based approximations of $\q$ can be erroneous (even with best-fit coefficients), Demo et al.\ \cite{demo2023deeponet} proposed a multi-fidelity modeling approach to learn the projection (truncation) error in POD reconstruction of steady-state flow problems. The leading POD coefficients are computed either using interpolation from pre-collected data or through sensor measurements. Then, a DeepONet is used to estimate the difference between the true field $\q$ and its reduced version $\qbar$. Geneva and Zabaras \cite{geneva2020a} proposed a multi-fidelity generative model that produces high-fidelity realizations of turbulent flows by conditioning an invertible neural network on the corresponding low-fidelity predictions.

Particularly relevant to the scope of this review is the closure learning methodology in \cite{ahmed2023a}, where $\G_{\pcal}(\qred;\pphys) = \F(\qred;\pphys)$ is considered as the low-fidelity model for the reduced quantity (i.e., ignoring the commutator error). Then, a multi-fidelity DeepONet is trained to learn the contribution of the closure term in a predictor-corrector fashion. The authors combined this multi-fidelity viewpoint with a posteriori learning to enhance the stability of the resulting predictions. The \emph{model fusion} approach in \cite{pawar2021model} can also be interpreted as a multi-fidelity learning approach, where the predictions of the Galerkin-POD model are used to guide a deep neural network toward an improved estimate of the leading POD coefficients. Sen et al.\ \cite{sen2018evaluation} considered a multi-fidelity modeling approach to address the closure problem in multiscale systems with distinct meso- and macroscales. They considered a shocked particulate flow where the mesoscale effects are incorporated into the macroscale model in the form of a homogenized drag value. A low-fidelity surrogate model for the drag is first developed by the modified Bayesian Kriging method using a set of low-fidelity mesoscale simulations. A multi-fidelity learning framework is subsequently constructed to correct the low-fidelity predictions needing only a few high-fidelity simulations. 

It is worth noting that using a multi-fidelity learning approach over a single fidelity one should be performed carefully. For example, \cite{giselle2019issues} showed that including low-fidelity samples along with the high-fidelity samples can lead to less accurate surrogates than just using the available high-fidelity samples.

\subsection{Multiscale }




Any review of closure modeling would not be complete without referencing the classical methods for reduction of systems exhibiting time-scale separation between variables that are needed to fully describe the state of the system. A comprehensive review on various methods to perform reduction for such systems can be found in \cite{givon2004}, while more recent accounts are in \cite{pavliotis2008, e2011}. The main classes of methods are invariant manifolds \cite{constantin1989integral}, averaging \cite{sanders2007averaging, gear2003projective,abdulle2012heterogeneous}, and white noise approximation \cite{majda2001}. We note that these methods can also be interpreted as limiting cases of the more general MZ formalism (see, e.g.,\ \cite{givon2004}). 

In recent years, there have appeared data-driven variants of the approaches above for systems exhibiting time-scale separation. The main topics of these approaches are the discovery of homogenized/effective equations for slow variables \cite{arbabi2020linking, callaham2023}, the discovery of collective variables that promote time-scale separation \cite{vlachas2022, dsilva2016data, crabtree2023gans, liu2023data}, as well as the estimation of closure terms to account for the effect of fast variables on slow variables \cite{mou2021data, san2022variational, koc2022verifiability, mou2023efficient, crabtree2023gans, lee2023learning, callaham2023}.

\section{Outlook}\label{sec:conclusions} 
The field of turbulence closure with machine learning has advanced rapidly in recent years, but major challenges remain \cite{spalart2023}. Here we list a few of these, including directions for future research. 

\subsection{Interpretability and generalizability}\label{sec:interpretability}
The machine-learned closure models discussed in this paper are generally represented by neural networks that do not have a clear sense of interpretability. However, interpretability is a key requirement for machine-learned closure models to find their way to engineering practice. Consequently, recent research efforts \cite{schmelzer2020,zanna2020,mojgani2022,jakhar2023} are being directed on providing interpretable, closed-form expressions for the closure model, e.g.,\ using techniques like sparse linear regression with a physics-informed library of candidate functions. Besides physical interpretability, this approach also promises to require fewer training data and have lower training costs \cite{zanna2020}.

Next to interpretability, generalizability is an important research area that is a prerequisite for trained closure models to be adopted in various situations. Training a model for certain parameters (e.g.,\ Reynolds number), geometries, boundary conditions, and initial conditions and using it beyond this parameter set could be seen as the ``holy grail'' in closure modeling research. However, this is currently out of reach in its generality. Similarly, training a model with training data that is associated to a certain grid, time step, and discretization scheme and reusing it under different conditions is a difficult but central goal in closure modeling that is drawing a lot of attention. A key issue for extrapolation in time is dealing with the chaotic nature of many physical models. Learning chaotic dynamics with machine learning has been investigated for several years \revthree{\cite{pathak2018,vlachas2018,fan2020,li2022,vlachas2022,schiff2024}}, but needs to further extend into the closure modeling community.

\subsection{Reduced model forms and including known physics}
It is generally acknowledged that including known physics  in data-driven methods is beneficial, reducing the amount of training data required and improving generalizability. In the closure modeling context of this paper, one would similarly argue that the closure form \eqref{eqn:closure_model_form_ddt} is more physics-informed than the physics-agnostic form \eqref{eqn:closure_model_form} and therefore has more potential to generalize and to train. However, the evidence that supports this claim is still rather thin (some examples are given in \cite{melchers2023}), and some studies are pointing in the reverse direction \cite{kochkov2021,stachenfeld2022a,bose2024}. Restricting the reduced models to have a closure form or a particular known symmetry was \textit{limiting} their performance in these cases. \RC{rev3_comment38}\revthree{Similarly, so-called foundation models have been recently developed in the weather forecasting community that are completely agnostic of the physics (although there is an implicit dependence on the physics through the training data that is generated from physics-based solvers) \cite{bodnar2024}. These foundation models are then fine-tuned for a specific task. It is still an open question whether such foundation models could play a role in a closure model context, or whether they could even make the entire concept of closure models redundant. Perhaps a `foundation closure model' (associated to a certain PDE) could be constructed, which is fine-tuned for a specific prediction task (boundary conditions, geometry, parameters, etc.). However, for a PDE like the Navier-Stokes equations, this seems at present an insurmountable task, given the effectively inexhaustible richness of possible flow configurations. In contrast to the `top-down' approach of foundation models, `bottom-up' approaches like the building-block model  presented in \cite{arranz2024} have been developed. These new approaches are still in the early phase of development, and which approach will stand the test of time will likely depend on a multitude of factors, including the availability of high-quality training data and access to computational resources.}

\subsection{Stability of hybrid models}
Scientific machine learning for closure leads to hybrid models that couple traditional solvers with neural networks. The often black-box nature of the neural network component can hinder the study of stability properties of the hybrid model, and, even worse, can amplify the risk of catastrophic instabilities. The problem of instabilities for coupled models has long been known (see, e.g.,\ \cite{keyes2013multiphysics} and references therein) and is not specific to the use of neural network closures. For neural network-based closures, possible overfitting of the training data can exacerbate the instability problem \cite{um2020}. Approaches like a posteriori learning presented in Section \ref{sec:solution_fitting} that rely on training the neural network so that the hybrid model accurately predicts a rolled-out trajectory of certain length have shown promise. However, thorough investigation of the stability properties of such approaches is still an active research topic.


\subsection{Network architectures}
Choosing the right network architecture and optimizing its parameters is often more of an art than a 
science. In this paper, little attention has been paid to the different network architectures used in the context of closure modeling because new architectures are appearing continuously and because it is still highly non-trivial which one to choose. Early machine learning-based closure modeling efforts focused on simple neural network layouts (fully connected multi-layer perceptrons) \cite{gamahara2017,maulik2019}, while later convolutional neural networks and autoencoders became popular approaches. Convolutional networks are useful for structured data and mimic, to a certain extent, finite difference operators, which is why they work well in the context of differential equations. Autoencoders are appropriate in the context of compression, for example in reduced order models. More recently, graph neural nets \cite{sanchez-gonzalez2020a,brandstetter2023} have gained attention because they work on unstructured data, and neural operators because of their grid independence \cite{lutjens2022,lu2021,li2021,cao2023lno,xu2024equivariant,bartolucci2023}. Only a few studies investigate the effect of the network architecture and the associated optimization process, e.g.,\ \cite{taghizadeh2021}.






\subsection{Benchmarking, test cases, and datasets}

Currently, many papers consider simplified PDE problems as test cases, such as the Burgers' equation, the Kuramoto-Sivashinsky equation, the Lorenz '96 equations, or the incompressible Navier-Stokes equations on periodic domains with low Reynolds numbers. A number of efforts have been initiated to build structured datasets that can be used to train machine learning methods \cite{kanov2015johns,rasp2020weatherbench,shen2022library,takamoto2023,towne2023database,yu2024climsim}. Because neural network design and training often require a lot of experience and trial and error, well-defined benchmarks, documented datasets, and relevant error metrics are of utmost importance. Extension to more complex physics is still an open topic---in fluid flows, closure models for compressible flows with shock waves, for irregular domains, for reactive flows, and for high Reynolds number flows are a sparsely explored but important territory \cite{shankar2023differentiable,sirignano2020dpm}. \RC{rev2_comment9}\revtwo{Part of benchmarking also means a fair assessment on whether large neural network models (and the associated training time and data requirements) can outperform classical (PDE-based) solvers in terms of accuracy per operation count or runtime. The choice of the reference solver that we compare against and try to beat is also a critical aspect. Researchers might be tempted to compare a neural-PDE solver against a classical PDE solver that is significantly more accurate and more expensive to justify the accuracy-cost trade-off of the former. However, in many cases, one can find another low-fidelity PDE solver that is as accurate as the neural-PDE solver but less-expensive. Thus, a fair comparison should involve a Pareto front analysis that surveys all solvers that are at our disposal.} \RC{rev1_comment3}\revone{In addition, many data-driven methods (for closures or otherwise) are limited by errors resulting from distribution shifting when deployed on test data that is significantly different from the training data. Learning useful `invariant' input features may assist with this but that is as yet a largely unsolved problem in machine learning \cite{zhao19invariant}.}



\section*{Acknowledgements}
The work of BS is partially supported by the project “Discretize first, reduce next” (with project number VI.Vidi.193.105) of the research programme NWO Talent Programme Vidi and the project “Unraveling Neural Networks with Structure-Preserving Computing” (with project number OCENW.GROOT.2019.044), both financed by the Dutch Research Council (NWO).

The work of PS is partially supported by the U.S. Department of Energy (DOE), Office of Science, Advanced Scientific Computing Research (ASCR) program
under the Uncertainty
Quantification for Multifidelity Operator Learning (MOLUcQ) project (Project No. 81739), under the Scalable, Efficient and Accelerated Causal Reasoning Operators, Graphs and Spikes for Earth and Embedded Systems (SEA-CROGS) project (Project No. 80278) and under DOE’s Scientific Discovery through Advanced Computing (SciDAC) program via a partnership on
Earth system model development between the Office of Biological and Environmental Research (BER) and ASCR
(Project No. 79699). Pacific Northwest National Laboratory is a multi-program national
laboratory operated for the U.S. Department of Energy by Battelle Memorial Institute
under Contract No. DE-AC05-76RL01830.

RM is partially supported by the U.S. Department of Energy, Office of Science, Advanced Scientific Computing Research program under the data-intensive scientific machine learning project (FOA-2493) and by the U.S. Department of Energy, Office of Science, Fusion Energy Sciences program under the DeepFusion Accelerator for Fusion Energy Sciences in Disruption Mitigations project (FOA-2905).

\appendix




\end{document}